\documentclass[preprintnumbers,amsmath,amssymb,aps]{revtex4-1}
\usepackage{graphicx}
\usepackage{dcolumn}
\usepackage{booktabs}
\usepackage{threeparttable}
\usepackage{longtable}
\usepackage{subfigure}
\usepackage{rotating}
\usepackage{multirow}
\usepackage{epstopdf}
\usepackage{tabularx}
\usepackage{amsmath}
\usepackage{amsfonts} 
\usepackage{mathrsfs} 
\usepackage{graphics}
\usepackage{array}

\usepackage{latexsym}
\usepackage{url}
\usepackage{lineno}
\usepackage{graphicx}
\usepackage{dcolumn}
\usepackage{booktabs}
\usepackage{threeparttable}
\usepackage{longtable}
\usepackage{subfigure}
\usepackage{rotating}
\usepackage{multirow}
\usepackage{epstopdf}
\usepackage{tabularx}
\usepackage{amsmath,amssymb}
\usepackage{amsmath,amscd}
\usepackage{amsfonts}
\usepackage{natbib}
\usepackage{numcompress}
\begin{document}

\title{Analytical solutions from integral transforms for transient fluid flow in naturally fractured porous media with and without boundary flux}

\author{Luis X.~Vivas-Cruz} \email{lvivas@posgrado.cidesi.edu.mx}
\affiliation{Centro de Ingenier\'ia y Desarrollo Industrial (CIDESI), Av. Playa Pie de la cuesta 702, Desarrollo San Pablo, Quer\'etaro, Qro 76125, Mexico}

\author{Jorge Adri\'an Perera-Burgos} \email{jorge.perera@cicy.mx}
\affiliation{CONACYT - Unidad de Ciencias del Agua, Centro de Investigaci\'on Cient\'ifica de Yucat\'an A.C. - Calle 8, No. 39, Mz. 29, S.M. 64, C.P. 77524, Canc\'un, Quintana Roo, M\'exico.}

\author{Alfredo Gonz\'alez-Calder\'on} \email{alfredo.gonzalez@cidesi.edu.mx}
\affiliation{CONACyT - CIDESI, Av. Playa Pie de la cuesta 702, Desarrollo San Pablo, Quer\'etaro, Qro 76125, Mexico}

\begin{abstract}
A kind of problems of radially symmetric transient fluid flow in a medium with a geometry similar to a hollow-disk can be addressed using the finite Hankel transform. However, the inverse Hankel transform [G.~Cinelli, Int.~J.~Engng.~Sci.,~{\bf 3}, 539 (1965)] works well only for homogeneous boundary conditions. We use the finite Hankel transform, together with the Laplace transform, to solve partial differential equations with inhomogeneous boundary conditions. With this aim, we propose a method to obtain an analytical solution of the problem of fluid flow in a finite naturally fractured reservoir with inner and outer boundaries having constant and time-dependent conditions, respectively. We assume that the reservoir has a producing well, with either constant terminal pressure or constant terminal rate, while the outer boundary has either an influx recharge or constant pressure. Using these case studies, we show that a part of the inverse transformation given by Cinelli can be expressed as closed formulas for long time solutions, at the same time that make it possible to capture the inhomogeneous condition and speed up the convergence of solutions. For Neumann-Neumann conditions, we show that the inverse expression of Cinelli is incomplete. In addition, an analysis of the flow characteristic curves in a reservoir with influx recharge is presented. The drawdown pressure curves are used to elucidate the statement that the pressure drop of a single-porosity reservoir with influx recharge resembles the flow behavior in a double-porosity closed reservoir, establishing a criteria to distinguish between both.

\noindent\\
{\bf Keywords}\\
Fluid flow in naturally fractured reservoir; Influx recharge at the outer boundary; Joint Laplace-Hankel transform.
\end{abstract}
\maketitle
\section{Introduction}
In groundwater science and petroleum engineering, the modelling of fluid flow in underground reservoirs has impact on project planning and reserve estimates. However, current models for fluid flow in reservoirs have limitations that affect their accuracy when they are applied in the tasks just mentioned. Therefore, there is a need for including in the governing equations the natural properties such as storage, porosities, permeabilities, wellbore storage, skin factor, or recharge. Furthermore, new mathematical developments with applications in pumping or well tests~(\cite{ozkan1988some,liu1990exact,chen1990analytical,young1992pressure,wu2002approximate,de2011analytical}) allows to understand phenomena, which could be challenging otherwise~(\cite{da1990well,singhal2010applied,nie2012dual,pedretti2014apparent,molinari2015analysis,
kuhlman2015multiporosity,zhou2019revisiting}).

Many fluid flow models have exact solutions in Laplace space, but their inverse transforms can be quite complex to obtain by means of contour integration in the complex plane~(\cite{yao2012transient,gonzalez2017exact}). Remarkably, the Hankel transform provides a simple way to treat radially symmetric problems, since their inverse transform formulas are the solutions of the models~(\cite{sneddon1946iii,cinelli1965extension,jiang2010general}). One of these solutions are given by~\citet{cinelli1965extension}; nevertheless, because their relationships are for homogeneous boundary conditions, they are strongly limited in their application to describe fluid flow in reservoirs associated with a hollow disk geometry. In order to extend the applicability of the finite Hankel transform derived by Cinelli, in Refs.~(\cite{xi1991theoretical}) and (\cite{wang1992elastodynamic}) developed a mathematical procedure that involves its use for inhomogeneous boundary conditions. In those studies the initial and boundary value problem are expressed as the sum of a dynamic part with inhomogeneous boundary conditions and a quasi-static part with homogeneous boundary conditions, in such a way that the solution from the quasi-static differential equation can make use of the mentioned formulas. In a similar fashion, in this work we divided the solution in a stationary and a transient part, which appropriately allows to solve the problem with inhomogeneous conditions. By contrast, in Refs.~(\cite{xi1991theoretical}) and (\cite{wang1992elastodynamic}) is solved a hyperbolic model for describing displacement in elastodynamics, while in this study we solve a parabolic model for describing the fluid flow in a double-porosity reservoir. Models of fluid flow in reservoirs with inhomogeneous boundary conditions, for example, influx recharge, have been the subject of different studies~(\cite{doublet1995decline,del2014pressure,wang2017transient}), but, to our best knowledge, they have not been analytically solved for double-porosity systems.

In addition to the Hankel transform, in order to solve partial differential equations, the Laplace transforms can be jointly used, taking us to the joint Laplace-Hankel transform or JLHT~(\cite{poularikas2010transforms,debnath2014integral}). An application in models of fluid flow in reservoirs is found in~\citet{babak2014unified}, where finite and infinite reservoirs are considered, each of them having a centered well with an infinitely small radius. For hollow-disk geometry, the finite Hankel transform was used to solve a triple-porosity fluid flow model, with a constant pressure and zero flux at the inner and outer boundaries, respectively, and considering a non-zero well radius~(\cite{clossman1975aquifer}). Also, the JLHT has been used in the study of crossflow in stratified systems; for instance,~\citet{boulton1977unsteady} provided the relationships of flow through horizontal layers of a fissured aquifer, restricted to have vertical permeability, and with a wellbore represented as a line source pumping at a constant rate. A similar system, with a partially penetrating well, is found in a study by~\citet{javandel1983analytical}. In the early 60s, \citet{katz1962theoretical} and \citet{russell1962performance} were the pioneers in the study of crossflow in stratified reservoirs using the Fourier and the Hankel transforms, respectively. Subsequently, \citet{prats1986interpretation} showed that stratified reservoir and single-layer reservoirs have similar behavior for large periods of time. More complex systems were analyzed by \citet{shah1992transient}, they included two flowing intervals in a partial completion well. On the other hand, exact solutions by means of other mathematical procedures can be consulted in Refs.~(\cite{gomes1993analytical,ehlig1987new}), for layered aquifers, and in (\cite{matthews1967pressure,chen1989transient}), for oil reservoirs. Additional applications of JLHT are found in~(\cite{carslaw1959conduction,poularikas2010transforms,debnath2014integral}).

There is a lack of exact solutions using the JLHT, for models of fluid flow in bounded or infinite reservoirs. Partly, this may be due to the complexities inherent to the method that will be exposed in this work. In order to contribute to the studies in this direction, we use the JLTH to solve the double-porosity model of~\citet{warren1963behavior} using the following combinations of specific boundary conditions (BCs): Dirichlet-Dirichlet (DD), Dirichlet-Neumann (DN), Neumann-Dirichlet (ND), and Neumann-Neumann (NN). The inner condition is given by either a constant terminal pressure or constant terminal rate. Meanwhile, the outer boundary has a constant pressure or has a flux defined as a ``Ramp" rate function to simulate natural water influx or slow-starting waterfloods from injector wells~(\cite{doublet1995decline}). Also, this latter function can be interpreted as rock heterogeneities at the outer boundary that obstructs the flow channels~(\cite{del2014pressure,doublet1995decline}). We note that our solutions generalize the relationships of fluid flow in a single-porosity medium, which were released in other works, and can be found in (\cite{muskat1934flow}), for DD-BCs; (\cite{muskat1934flow,hurst1934unsteady}), for DN-BCs; (\cite{matthews1967pressure}), for ND-BCs; and in (\cite{muskat1934flow,matthews1967pressure,del2014pressure}), for NN-BCs.

Warren and Root model has been widely used in well and pumping tests analy\-sis~(\cite{bourdet1989use,kruseman1994analysis,gringarten2008straight,ahmed2011advanced}). Nevertheless, the NN-BCs case deserves a special attention, due to the lack of studies in this regard, in such a way that it is important to provide type curves that include the effects of influx recharge at the outer boundary. Indeed, it has been observed that a single-porosity reservoir with influx recharge has a characteristic behavior similar to that of a double-porosity reservoir without influx recharge~(\cite{doublet1995decline,del2014pressure}). By extension, it could be expected that a double-porosity reservoir with influx recharge has a behavior similar to that given by a triple-porosity reservoir without influx recharge. This observation is attended in this work in order to elucidate whether the drawdown pressures of models with recharge and without it can be considered equivalent.

The contribution of this work is threefold:
\begin{enumerate}
\item Obtain analytical solutions for the aforementioned study cases.
\item Present characteristic behaviors of the drawdown pressure and flux for the study cases.
\item Show the similarities and differences between a model with influx recharge and without it, but this latter case with an additional porosity.
\end{enumerate}

This work is organized as follows: In Sections~\ref{Section_Mathematical_model} and \ref{Section_boundary conditions}, the flow model and the boundary conditions are presented, respectively. In Section~\ref{Section_Hankel-Laplace_Sol}, the procedure to find the exact solutions is given. Subsequently, in Section~\ref{Section_discussion}, a numerical validation using the Stehfest method is carried out, the convergence of the Cinelli solutions to the exact result is analyzed, and a discussion about the stationary solutions is also presented. In Section~\ref{Section_characteristic_behavior_of_solutions}, the characteristic curves of the fluid flow model are exposed, and a comparison is done between models with influx recharge and without it. Finally, in Section~\ref{Section_conclusions}, some general conclusions are drawn.

\section{Flow model}\label{Section_Mathematical_model}

\begin{figure}[ht]
{\includegraphics[width=0.85\linewidth]{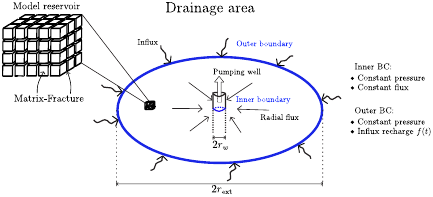}}
\caption{Schematic representation of the reservoir and the fluid flow. Drainage area of reservoir is illustrated by a big circle of radius $r_{\text{ext}}$, while in its center there is a cylinder of radius $r_{\text{w}}$ that represent the pumping well. Other boundary conditions used in this work are specified on the right of the figure
}
\label{ReservoirModel}
\end{figure}

Radial fluid flow in a double-porosity medium is described by the model of~\citet{warren1963behavior}. This model considers a slightly compressible fluid through two overlapping porous systems, matrix and fractures. The matrix discharges into the fractures, and the fractures carry the flow toward the wellbore. The matrix has a low permeability and a high storage capacity, while the fracture system has a high permeability and a low storage capacity. It is also assumed a homogeneous and isotropic porous system. An illustration of the reservoir and its fluid flow is given in Fig.~\ref{ReservoirModel}.

Warren and Root model considers an equation for the pressure $h_2$ in fractures,
\begin{linenomath*}
\begin{equation}\label{Eq_fractures}
\frac{k_2}{\mu} \frac{1}{r}\frac{\partial}{\partial r}\left(r \frac{\partial h_2}{\partial r}\right) -\phi_1 c_1 \frac{\partial h_1}{\partial t} = \phi_2 c_2 \frac{\partial h_2}{\partial t}, 
\end{equation}
\end{linenomath*}
and an equation for the pressure $h_1$ in matrix blocks,
\begin{linenomath*}
\begin{equation}\label{Eq_matrix}
\phi_1 c_1 \frac{\partial h_1}{\partial t} = \frac{\alpha \kappa_1}{\mu}(h_2 - h_1).
\end{equation}
\end{linenomath*}
In the previous equations, subscripts 1 and 2 indicate the matrix medium and the fracture medium, respectively; $\kappa_i$, $\phi_i$ and $c_i$, are the permeabilities, porosities and total compressibilities of medium $i$; $\mu$ is the fluid viscosity; and $\alpha$ is the shape factor.

Eqs.~(\ref{Eq_fractures}) and (\ref{Eq_matrix}) in reduced units are as follows:
\begin{linenomath*}
\begin{eqnarray}\label{2Pmodel}
\displaystyle \omega \frac{\partial h_{2D}}{\partial t_{D}} + (1-\omega)\frac{\partial h_{1D}}{\partial t_{D}}  &=&  \displaystyle \frac{\partial^2 h_{2D}}{\partial r_D^2} + \frac{1}{r_D}\frac{\partial h_{2D}}{\partial r_D}, \qquad 1\leq r_D \leq r_{D\text{ext}},\\
\displaystyle (1-\omega)\frac{\partial h_{1D}}{\partial t_{D}} &=&  \lambda(h_{2D} - h_{1D}),
\end{eqnarray}
\end{linenomath*}
where the dimensionless dependent variables are given by
\begin{linenomath*}
\begin{equation}\label{Dim_variables_PressureandFlux}
\displaystyle h_{iD} \equiv \left\lbrace \begin{array}{lll}
\displaystyle \frac{h_0 - h_i(r,t)}{h_0 - h_w}, && \mbox{for constant pressure},\\[15pt]
\displaystyle \frac{2\pi H \kappa_2}{\mu q}[h_0 - h_i(r,t)], && \mbox{for constant flux}.
\end{array}\right.
\end{equation}
\end{linenomath*}
The dimensionless independent variables are defined as
\begin{linenomath*}
\begin{equation}\label{dimensionless_variables_1}
\begin{split}
t_D \equiv \frac{\kappa_2 t}{\mu r_w^2(\phi_1 c_1 + \phi_2 c_2)}, \quad r_D \equiv \frac{r}{r_w},
\end{split}
\end{equation}
\end{linenomath*}
and the parameters of the model are
\begin{linenomath*}
\begin{equation}\label{dimensionless_variables_2}
\begin{split}
\omega \equiv \frac{\phi_2 c_2}{\phi_1 c_1 + \phi_2 c_2}, \quad \lambda \equiv \frac{\alpha r_w^2 \kappa_1}{\kappa_2}, \quad r_{D\text{ext}} \equiv \frac{r_{\text{ext}}}{r_w}.
\end{split}
\end{equation}
\end{linenomath*}
In these latter definitions, $r_w$ is the well radius; $r_{D\text{ext}}$ is the dimensionless outer radius; $\omega$ is the fracture storage coefficient; $\lambda$ is the interporosity flow coefficient; $h_0$ and $h_w$ are the reference pressure and the pressure at the bottomhole, respectively; $H$ is the thickness of the uniform horizontal formation; and $q$ is the constant volumetric flow rate.

In addition, the mass flux $j_2$ in reduced units is defined as follows~(\cite{matthews1967pressure}):
\begin{linenomath*}
\begin{equation}\label{flux_Darcylaw}
j_{2D}(t_D) = -r_D\frac{\partial h_{2D}(r_D,t_D)}{\partial r_D}\Bigr|_{r_D=1},
\end{equation}
\end{linenomath*}
where
\begin{linenomath*}
\begin{equation}\label{Flux_dimensionless}
\displaystyle j_{2D} = \left\lbrace \begin{array}{lll}
\displaystyle \frac{\mu r_w}{\rho \kappa_2 (h_w-h_0)}j_2, && \mbox{for constant pressure},\\[15pt]
\displaystyle -\frac{2\pi H \kappa_2}{\mu q}j_2, && \mbox{for constant flux},
\end{array}\right.
\end{equation}
\end{linenomath*}
where $\rho$ is the density of the fluid in the fractures per unit volume.
\section{Initial and boundary conditions}\label{Section_boundary conditions}
Assuming that the reservoir has a constant pressure at time zero, the initial conditions are written as
\begin{linenomath*}
\begin{equation}\label{initial_conditions}
h_{1D}(r_D,0) = h_{2D}(r_D,0) = 0.
\end{equation}
\end{linenomath*}
The boundary condition at the bottomhole, when a constant pressure is imposed, is given by
\begin{linenomath*}
\begin{equation}\label{constant_head}
h_{2D}(1,t_D) = 1.
\end{equation}
\end{linenomath*}
On the other hand, when a constant flow at the bottomhole is imposed and an influx recharging the reservoir through the outer boundary is considered, the boundary conditions are as follows:
\begin{linenomath*}
\begin{equation}\label{constant_rate}
\displaystyle r_D \frac{\partial h_{2D}(r_D,t_D)}{\partial r_D} = \left\lbrace \begin{array}{lll}
\displaystyle -1, && \mbox{for $r_D=1$,}\\[15pt]
\displaystyle f(t_D) = -q_{\text{ext}}(1-\mbox{e}^{-t_D/\gamma}), && \mbox{for $r_D=r_{D\text{ext}}$},
\end{array}\right.
\end{equation}
\end{linenomath*}
where $q_{\text{ext}}\geq 0$ is the influx factor, and $\gamma$ is a parameter to change the slope of the ``Ramp" rate function~(\cite{doublet1995decline,del2014pressure}). Note that $q_{\text{ext}}=0$ is for a reservoir with zero flux at the outer boundary~(\cite{van1949application}). 

A constant pressure at the outer boundary is also considered:
\begin{linenomath*}
\begin{equation}\label{rext_boundary}
h_{2D}(r_{D\text{ext}},t_D) = 0.
\end{equation}
\end{linenomath*}
Some combinations of these BCs are shown in Table \ref{TableConditions}. They comprise the case studies analyzed in this work.

\renewcommand{\arraystretch}{0}
\begin{table}[t]
\small\addtolength{\tabcolsep}{-10pt}
\caption{\label{TableConditions}Case studies and their dimensionless boundary conditions used to solve the model~(\ref{2Pmodel}).}
\begin{tabular}{p{1.5cm}p{4cm}p{6.5cm}p{6.5cm}}
\toprule[0.7pt]
\begin{center}Case\end{center} & \begin{center}Boundary conditions\end{center} & \begin{center}Inner boundary\end{center} & \begin{center}Outer boundary\end{center}\\[0.5pt]
\cmidrule[0.7pt](r){1-4} \\[5pt]
\begin{center}
DD-BCs
\end{center}

&
\begin{center}
Dirichlet-Dirichlet
\end{center}
&
\begin{equation*}
\displaystyle h_{2D}(1,t_D) = 1
\end{equation*}

&
\begin{equation*}
\displaystyle h_{2D}(r_{D\text{ext}},t_D) = 0
\end{equation*}

\\[5pt]
\begin{center}
DN-BCs
\end{center}

&
\begin{center}
Dirichlet-Neumann
\end{center}
&
\begin{equation*}
\displaystyle h_{2D}(1,t_D) = 1
\end{equation*}

&
\begin{equation*}
\displaystyle r_D\frac{\partial h_{2D}(r_D,t_D)}{\partial r_D}\Bigr|_{r_D=r_{D\text{ext}}} = f(t_D)
\end{equation*}
\\[5pt]

\begin{center}
ND-BCs
\end{center}

&
\begin{center}
Neumann-Dirichlet
\end{center}
&
\begin{equation*}
\displaystyle r_D \frac{\partial h_{2D}(r_D,t_D)}{\partial r_D}\Bigr|_{r_D=1} = -1
\end{equation*}

&
\begin{equation*}
\displaystyle h_{2D}(r_{D\text{ext}},t_D) = 0
\end{equation*}
\\[5pt]
\begin{center}
NN-BCs
\end{center}

&
\begin{center}
Neumann-Neumann
\end{center}
&
\begin{equation*}
\displaystyle r_D \frac{\partial h_{2D}(r_D,t_D)}{\partial r_D}\Bigr|_{r_D=1} = -1
\end{equation*}
&
\begin{equation*}
\displaystyle r_D \frac{\partial h_{2D}(r_D,t_D)}{\partial r_D}\Bigr|_{r_D=r_{D\text{ext}}} = f(t_D)
\end{equation*}
\\
\bottomrule[0.7pt]
\end{tabular}
\end{table}

\section{Integral transform solutions}\label{Section_Hankel-Laplace_Sol}
In this section, using the JLHT and their inversion formulas, solutions of fluid flow in a double-porosity medium, Eqs.~(\ref{2Pmodel}), for the initial-boundary conditions in Section~\ref{Section_boundary conditions}, are presented. To validate the solutions, comparisons with results from Stehfest method are made in Section~\ref{Section_discussion}. Details of the calculations can be found in Appendix A. Henceforth, for simplicity in notation, we omit the subscript $D$ of the dimensionless variables previously defined.

Eqs.~(\ref{2Pmodel}) in Laplace space are written as
\begin{linenomath*}
\begin{equation}\label{2Pmodel_Laplace}
\begin{split}
\displaystyle \omega s \widehat{h}_{2} + (1-\omega)s \widehat{h}_{1}  & =  \displaystyle \frac{\partial^2 \widehat{h}_{2}}{\partial r^2} + \frac{1}{r}\frac{\partial \widehat{h}_{2}}{\partial r}, \qquad 1\leq r \leq r_{\text{ext}},\\
\displaystyle (1-\omega)s \widehat{h}_{1} & =  \lambda(\widehat{h}_{2} - \widehat{h}_{1}),
\end{split}
\end{equation}
\end{linenomath*}
and the BCs in Table \ref{TableConditions} are summarized as follows:
\begin{linenomath*}
\begin{eqnarray}
\displaystyle \widehat{h}_{2}(1,s) & = & \frac{1}{s}\label{condition1},\\
\displaystyle r \frac{\partial \widehat{h}_{2}(r,s)}{\partial r}\Bigr|_{r=1} & = & -\displaystyle \frac{1}{s}\label{condition2},\\
\displaystyle \widehat{h}_{2}(r_{\text{ext}},s) & = & 0\label{condition3},\\
\displaystyle r\frac{\partial \widehat{h}_{2}(r,s)}{\partial r}\Bigr|_{r=r_{\text{ext}}} & = & \widehat{f}(s)\label{condition4},\\
\displaystyle \widehat{h}_{i}(r,s) & = & 0, \qquad i=1,2\label{condition5},
\end{eqnarray}
\end{linenomath*}
where $s$ is the Laplace transform variable, $\widehat{x}$~denotes the Laplace transform of $x$, and
\begin{linenomath*}
\begin{eqnarray}\label{f(t)_Laplace}
\widehat{f}(s) & = & -q_{\text{ext}} \left( \frac{1}{s} - \frac{1}{s+1/\gamma}\right).
\end{eqnarray}
\end{linenomath*}
From Eqs.~(\ref{2Pmodel_Laplace}), the Laplace transform of pressure in the fractures obeys the modified Bessel ordinary differential equation:
\begin{linenomath*}
\begin{equation}\label{Laplace_2Dmodel_Bessel}
\frac{d^2 \widehat{h}_{2}}{d r^2} + \frac{1}{r}\frac{d \widehat{h}_{2}}{d r} - \eta(s) \; \widehat{h}_{2} = 0,
\end{equation}
\end{linenomath*}
where
\begin{linenomath*}
\begin{equation}\label{eta(s)}
\eta(s) = \frac{s\omega(1-\omega) +\lambda}{s(1-\omega) + \lambda}s.
\end{equation}
\end{linenomath*}

We use the finite Hankel transform to obtain an analytical solution of our study model, so by taking this transform to Eq.~(\ref{Laplace_2Dmodel_Bessel}) leads to the following expressions in the joint Laplace-Hankel space:
\begin{linenomath*}
\begin{subequations}\label{h_2(k,s)}
\begin{align}
\displaystyle \tilde{\widehat{h}}_{2}(k_i,s) =& \displaystyle \mathcal{F}(k_i,s) \left(\widehat{h}_{2}(r,s)\Bigr|_{r=r_{\text{ext}}} \frac{J_0(rk_i)|_{r=1}}{J_0(rk_i)|_{r=r_{\text{ext}}}} - \widehat{h}_{2}(r,s)\Bigr|_{r=1}\right), \quad \mbox{for DD-BCs}, \label{h_2(k,s)a}\\
=& \displaystyle \mathcal{F}(k_i,s) \left( r\frac{\partial \widehat{h}_{2}}{\partial r}\Bigr|_{r=r_{\text{ext}}} \frac{J_0(rk_i)|_{r=1}}{k_i J'_0(rk_i)|_{r=r_{\text{ext}}}} - \widehat{h}_{2}(r,s)\Bigr|_{r=1}\right), \quad \mbox{for DN-BCs},\label{h_2(k,s)b}\\
=& \displaystyle \mathcal{F}(k_i,s) \left(\widehat{h}_{2}(r,s)\Bigr|_{r_{\text{ext}}} \frac{J'_0(rk_i)|_{r=1}}{J_0(rk_i)|_{r_{\text{ext}}}} - \frac{1}{k_i}r\frac{\partial \widehat{h}_{2}}{\partial r}\Bigr|_{r=1}\right), \quad \mbox{for ND-BCs},\label{h_2(k,s)c}\\
=& \displaystyle \mathcal{F}(k_i,s) \left( r\frac{\partial \widehat{h}_{2}}{\partial r}\Bigr|_{r_{\text{ext}}} \frac{J'_0(rk_i)|_{r=1}}{k_i J'_0(rk_i)|_{r_{\text{ext}}}} - \frac{1}{k_i}r\frac{\partial \widehat{h}_{2}}{\partial r}\Bigr|_{r=1}\right), \quad \mbox{for NN-BCs},\label{h_2(k,s)d}
\end{align}
\end{subequations}
\end{linenomath*}
where $\tilde{x}$ is the finite Hankel transform of $x$~\cite{cinelli1965extension} and $\mathcal{F}(k_i,s) = \displaystyle 2/\{ \pi [\eta(s) + k_i^2] \}$. After that, we take the inverse Laplace transform of Eqs.~(\ref{h_2(k,s)}) and then we use the inverse finite Hankel transform~(\cite{cinelli1965extension}) to get:

\begin{linenomath*}
\begin{subequations}\label{inverse_h_2(k,s)}
\begin{align}
\displaystyle h_{2}(r,t) =& \displaystyle \frac{\pi^2}{2} \sum_{k_i>0}\frac{k_i^2 J_0^2(k_i r_{\text{ext}})\tilde{h}_{2}(k_i,t)}{J_0^2(k_i)-J_0^2(k_ir_{\text{ext}})}\mathcal{I}_{0,0}(k_i,r,1), \quad \mbox{for DD-BCs}, \label{inverse_h_2(k,s)a} \\
=& \displaystyle \frac{\pi^2}{2} \sum_{k_i>0}\frac{k_i^2 J_1^2(k_i r_{\text{ext}})\tilde{h}_{2}(k_i,t)}{J_0^2(k_i)-J_1^2(k_ir_{\text{ext}})}\mathcal{I}_{0,0}(k_i,r,1), \quad \mbox{for DN-BCs}, \label{inverse_h_2(k,s)b} \\
=& \displaystyle \frac{\pi^2}{2} \sum_{k_i>0}\frac{k_i^2 J_0^2(k_i r_{\text{ext}})\tilde{h}_{2}(k_i,t)}{J_1^2(k_i)-J_0^2(k_ir_{\text{ext}})}\mathcal{I}_{1,0}(k_i,1,r),\quad \mbox{for ND-BCs}, \label{inverse_h_2(k,s)c} \\
=& \displaystyle \frac{\pi^2}{2} \sum_{k_i>0}\frac{k_i^2 J_1^2(k_i r_{\text{ext}})\tilde{h}_{2}(k_i,t)}{J_1^2(k_i)-J_1^2(k_ir_{\text{ext}})}\mathcal{I}_{1,0}(k_i,1,r),\quad \mbox{for NN-BCs}, \label{inverse_h_2(k,s)d}
\end{align}
\end{subequations}
\end{linenomath*}
where $\mathcal{I}_{m,n}(\text{x},\text{y},\text{z}) = J_m(\text{x}\text{y})Y_n(\text{x}\text{z}) - Y_m(\text{x}\text{y})J_n(\text{x}\text{z})$ and $k_i$ are the roots of $\mathcal{I}_{0,0}(k_i,1,r_{\text{ext}})$, $\mathcal{I}_{1,0}(k_i,r_{\text{ext}},1)$, $\mathcal{I}_{1,0}(k_i,1,r_{\text{ext}})$, and $\mathcal{I}_{1,1}(k_i,1,r_{\text{ext}})$, for Eq.~(\ref{inverse_h_2(k,s)a}),~(\ref{inverse_h_2(k,s)b}),~(\ref{inverse_h_2(k,s)c}), and,~(\ref{inverse_h_2(k,s)d}), respectively. Eqs.~(\ref{inverse_h_2(k,s)}) are for homogeneous inner BCs, since $\mathcal{I}_{0,0}(k_i,1,r)=0$ and $\partial \mathcal{I}_{1,0}(k_i,1,r)/\partial r=0$ when $r=1$. Therefore, they require corrections when inhomogeneous BCs are considered. In Ref.~(\cite{xi1991theoretical}) a mathematical procedure is given to help with the correct application of the formulas of~\citet{cinelli1965extension}, when a hyperbolic differential equation is considered. However, in this work we choose to follow an alternative procedure, which only works when a long-time asymptotic (stationary) solution exists. Below a discussion and evidences are presented. We propose the following procedure to correct the Eqs.~(\ref{inverse_h_2(k,s)}):
 
\begin{itemize}
\item Regarding $\tilde{h}_{2}(k_i,t)$ in Eqs.~(\ref{inverse_h_2(k,s)}), this function is rewritten as follows
\begin{linenomath*}
\begin{equation}
\tilde{h}_{2}(k_i,t) = \tilde{g}(k_i,t) + \tilde{\phi}(k_i),
\end{equation}
\end{linenomath*}
where $\tilde{\phi}(k_i)$ contains the time-independent (stationary) terms. Therefore,  $\tilde{g}(k_i,t)$ is the transient part of the solution. For more details about $\tilde{g}(k_i,t)$ for each BCs case, see Appendix A.
 \item We write each equation in (\ref{inverse_h_2(k,s)}) as
\begin{linenomath*}
\begin{equation}\label{sum+quasiStationary}
h_{2}(r,t) = \mathcal{H}^{-1}\left[\sum \tilde{g}(k_i,t)\right](r) + \mathcal{H}^{-1}\left[\sum \phi(k_i)\right](r),
\end{equation}
\end{linenomath*}
where $\mathcal{H}^{-1}[\sum \tilde{g}(k_i,t)](r)$ has the same form that the infinite sum in Eqs.~(\ref{inverse_h_2(k,s)}), regarding that $\tilde{h}_{2}(k_i,t)$ is replaced by $\tilde{g}(k_i,t)$. This step is similar to the procedure in~\citet{xi1991theoretical}, where the solution is divided in two parts.
\item Since $\mathcal{H}^{-1}\left[\sum \phi(k_i)\right](r)$ is time independent, these terms are the stationary solution of the study model (\ref{2Pmodel}), i.e.~the solution to the Laplace ordinary differential equation,
\begin{linenomath*}
\begin{equation}\label{LaplaceODE}
\frac{1}{r}\frac{\partial}{\partial r}\left( r \frac{\partial h_{2,s}}{\partial r} \right)= 0,
\end{equation}
\end{linenomath*}
with the inhomogeneous BCs in Table \ref{TableConditions} at limit $t\rightarrow\infty$. Therefore,
\begin{linenomath*}
\begin{equation}
h_{2,s}(r) = \lim_{t\rightarrow\infty} h_2(r,t) = \mathcal{H}^{-1}\left[\sum \phi(k_i)\right](r).
\end{equation}
\end{linenomath*}
In this respect, it is worth noting that: 1) Our procedure can be applied because outer BC allows to obtain a stationary solution for DD-BCs, DN-BCs, and ND-BCs cases. 2) Because in specific conditions there is no solution of Eq.~(\ref{LaplaceODE}) for NN-BCs case, another alternative procedure to obtain the long-time solutions can be used. Namely, the long-time solution is obtained by expanding in series, about $s=0$, the equations in Table \ref{TableAppendix}, and then taking the inverse Laplace transformation of this results, which lead us to the desired solution. We denoted this solution as $h_{2,q}=h_{2,q}(r,t)$.
 \item In view of the discussion above, the solution of model (\ref{2Pmodel}) is
 \begin{linenomath*}
\begin{equation}\label{eqReplaceClosed}
h_{2}(r,t) = h_{2,s}(r) + \mathcal{H}^{-1}\left[\sum \tilde{g}(k_i,t)\right](r),
\end{equation}
\end{linenomath*}
for DD-BCs, DN-BCs and ND-BCs cases. Meanwhile, for NN-BCs case, the solution is
\begin{linenomath*}
\begin{equation}
h_{2}(r,t) = h_{2,q}(r,t) + \mathcal{H}^{-1}\left[\sum \tilde{g}(k_i,t)\right](r).
\end{equation}
\end{linenomath*}
\end{itemize}

Using the previous procedure (for additional details see Appendix A) the exact analytical solutions are:\\

\noindent
\emph{{\bf DD-BCs case}}

\begin{linenomath*}
\begin{equation}\label{solWR_case1}
h_{2}(r,t) = 1 - \frac{\log(r)}{\log(r_{\text{ext}})} + \frac{\pi}{2} \sum_{i=1}^{\infty} \frac{ \tilde{g}(k_i,t) \mathcal{I}_{0,0}(k_i,r,1) J_0^2(r_{\text{ext}}k_i)}{J_0^2(k_i) - J_0^2(r_{\text{ext}}k_i)}
\end{equation}
\end{linenomath*}
and
\begin{linenomath*}
\begin{equation}\label{solFlux_case1}
j_{2}(t) = \frac{1}{\log(r_{\text{ext}})} - \frac{\pi}{2} \sum_{i=1}^{\infty} \frac{k_i \tilde{g}(k_i,t) \mathcal{I}_{0,1}(k_i,1,1) J_0^2(r_{\text{ext}}k_i)}{J_0^2(k_i) - J_0^2(r_{\text{ext}}k_i)}.
\end{equation}
\end{linenomath*}

\noindent
\emph{{\bf DN-BCs case}}

\begin{linenomath*}
\begin{equation}\label{solWR_case2}
h_{2}(r,t) = 1 - q_{\text{ext}}\log(r) + \frac{\pi^2}{2} \sum_{i=1}^{\infty} \frac{ k_i^2 \tilde{g}(k_i,t) \mathcal{I}_{0,0}(k_i,r,1) J_1^2(r_{\text{ext}}k_i)}{J_0^2(k_i) - J_1^2(r_{\text{ext}}k_i)}
\end{equation}
\end{linenomath*}
and
\begin{linenomath*}
\begin{equation}\label{Radialflux_case2}
j_{2}(t) = q_{\text{ext}} - \frac{\pi^2}{2} \sum_{i=1}^{\infty} \frac{ k_i^3 \tilde{g}(k_i,t) \mathcal{I}_{0,1}(k_i,1,1) J_1^2(r_{\text{ext}}k_i)}{J_0^2(k_i) - J_1^2(r_{\text{ext}}k_i)}.
\end{equation}
\end{linenomath*}
\noindent
\emph{{\bf ND-BCs case}}

\begin{linenomath*}
\begin{equation}\label{solWR_case3}
h_{2}(r,t) = \log\left( \frac{r_{\text{ext}}}{r} \right) + \frac{\pi}{2} \sum_{i=1}^{\infty} \frac{ \tilde{g}(k_i,t) \mathcal{I}_{1,0}(k_i,1,r) J_0^2(r_{\text{ext}}k_i)}{k_i[J_1^2(k_i) - J_0^2(r_{\text{ext}}k_i)]}.
\end{equation}
\end{linenomath*}

\noindent
\emph{{\bf NN-BCs case}}

\begin{linenomath*}
\begin{equation}\label{solWR_case4}
\begin{split}
h_{2}(r,t) =&\frac{\pi^2}{2} \sum_{i=1}^{\infty} \frac{\tilde{g}(k_i,t) \mathcal{I}_{1,0}(k_i,1,r) J_1^2(r_{\text{ext}}k_i)}{J_1^2(k_i) - J_1^2(r_{\text{ext}}k_i)}+ \frac{2}{r_{\text{ext}}^2-1}\left( \frac{r^2}{4}+t \right)\\[3pt]
& - \frac{r_{\text{ext}}^2}{r_{\text{ext}}^2-1}\log (r) - \frac{3r_{\text{ext}}^4 - 4r_{\text{ext}}^4 \log(r_{\text{ext}}) - 2r_{\text{ext}}^2 - 1}{4(r_{\text{ext}}^2-1)^2} + q_{\text{ext}} \Big[ \frac{2}{r_{\text{ext}}^2-1}\left( \frac{r^2}{4} + t \right) - \frac{\log (r)}{r_{\text{ext}}^2-1}\\[3pt]
& - \frac{r_{\text{ext}}^4 + 2r_{\text{ext}}^2 - 4r_{\text{ext}}^2 \log (r_{\text{ext}}) -3}{4(r_{\text{ext}}^2-1)^2} \Big] + \frac{2 \gamma q_{\text{ext}}(1 - \mbox{e}^{-t/\gamma})}{r_{\text{ext}}^2-1}.
\end{split}
\end{equation}
\end{linenomath*}

\section{Analysis of the solutions}\label{Section_discussion}

In this section, the solutions are numerically validated and we show that Eqs.~(\ref{inverse_h_2(k,s)a}), (\ref{inverse_h_2(k,s)b}), and (\ref{inverse_h_2(k,s)c}), lead to exact results, except for the inner boundary where the homogeneous BCs are hold. While for NN-BCs case, we can see that Eq.~(\ref{inverse_h_2(k,s)d}) need correction terms to get the exact solution. On the other hand, solutions making use of closed relationships, Eqs.~(\ref{solWR_case1}), (\ref{solWR_case2}), (\ref{solWR_case3}), and (\ref{solWR_case4}), hold the inhomogeneous BCs given in Table~\ref{TableConditions} and their results match those from Stehfest method.

\subsection{Validation of the exact solutions using Stehfest method}
Next, the derived exact-analytic solutions of model (\ref{2Pmodel}) are numerically validated. With this goal the Stehfest method~(\cite{stehfest1970algorithm}) is used to take the inverse Laplace transform of equations in Table~\ref{TableAppendix}.

\begin{figure}[ht]
{\includegraphics[width=0.45\linewidth]{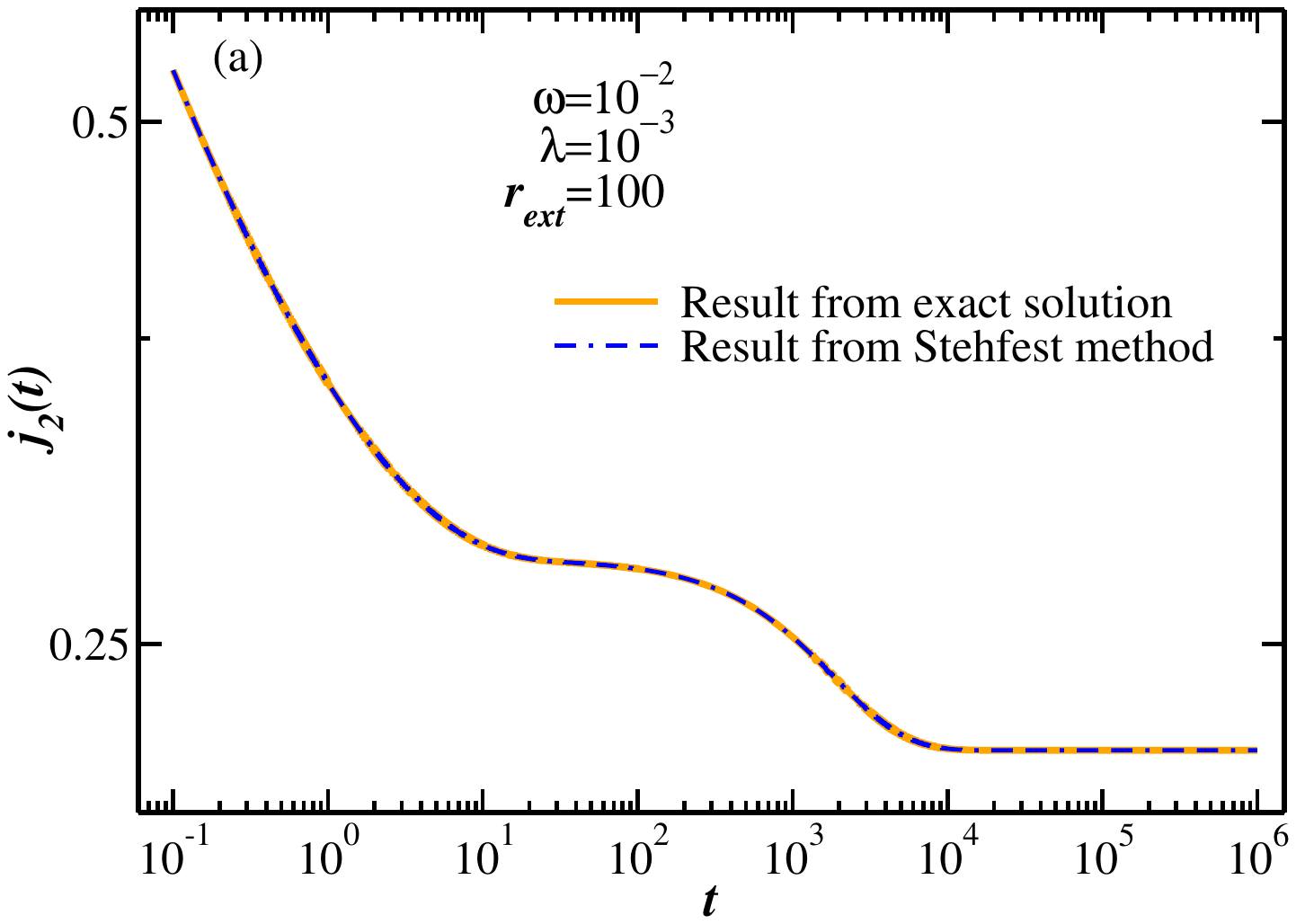}}
{\includegraphics[width=0.45\linewidth]{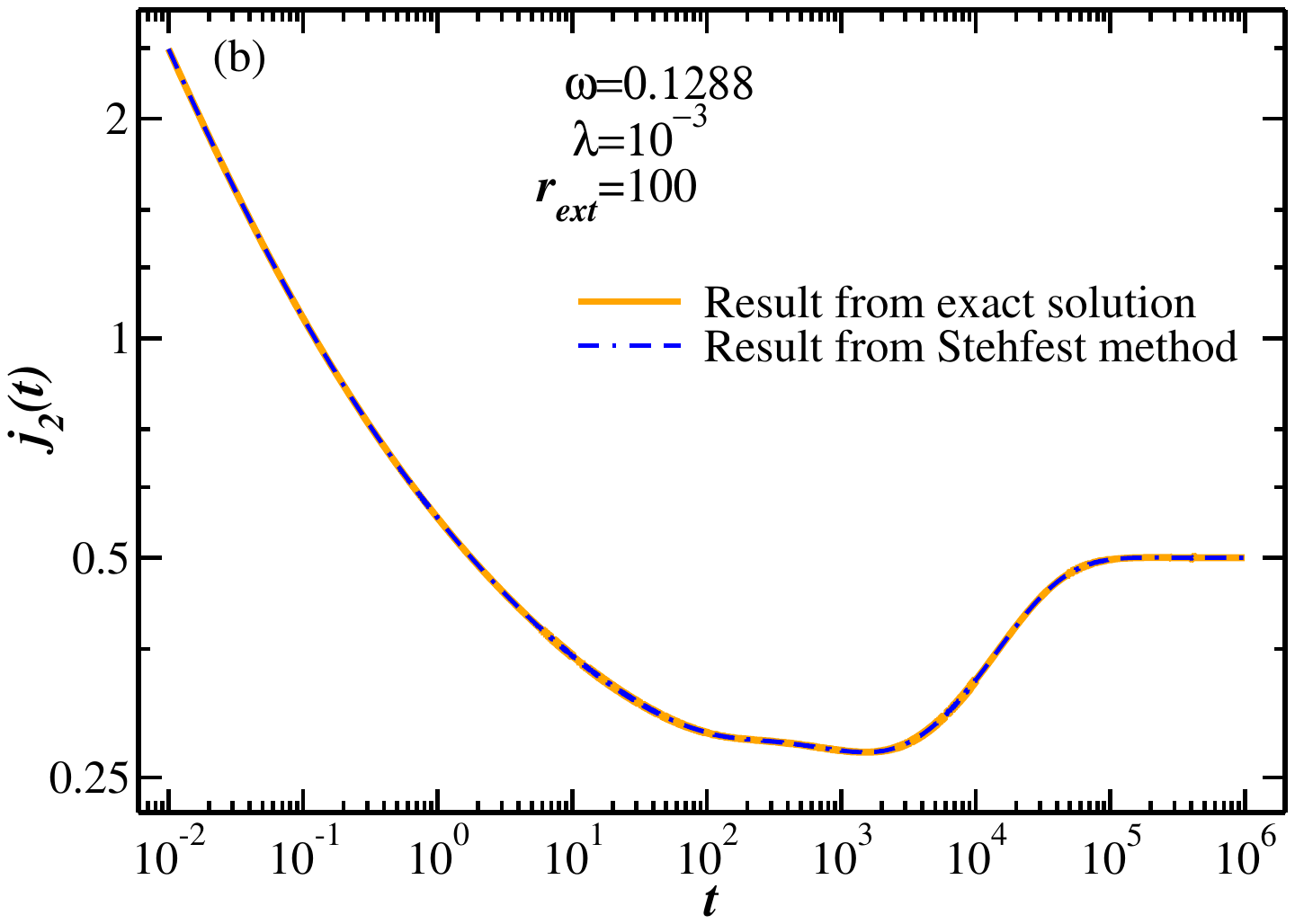}}\\
{\includegraphics[width=0.45\linewidth]{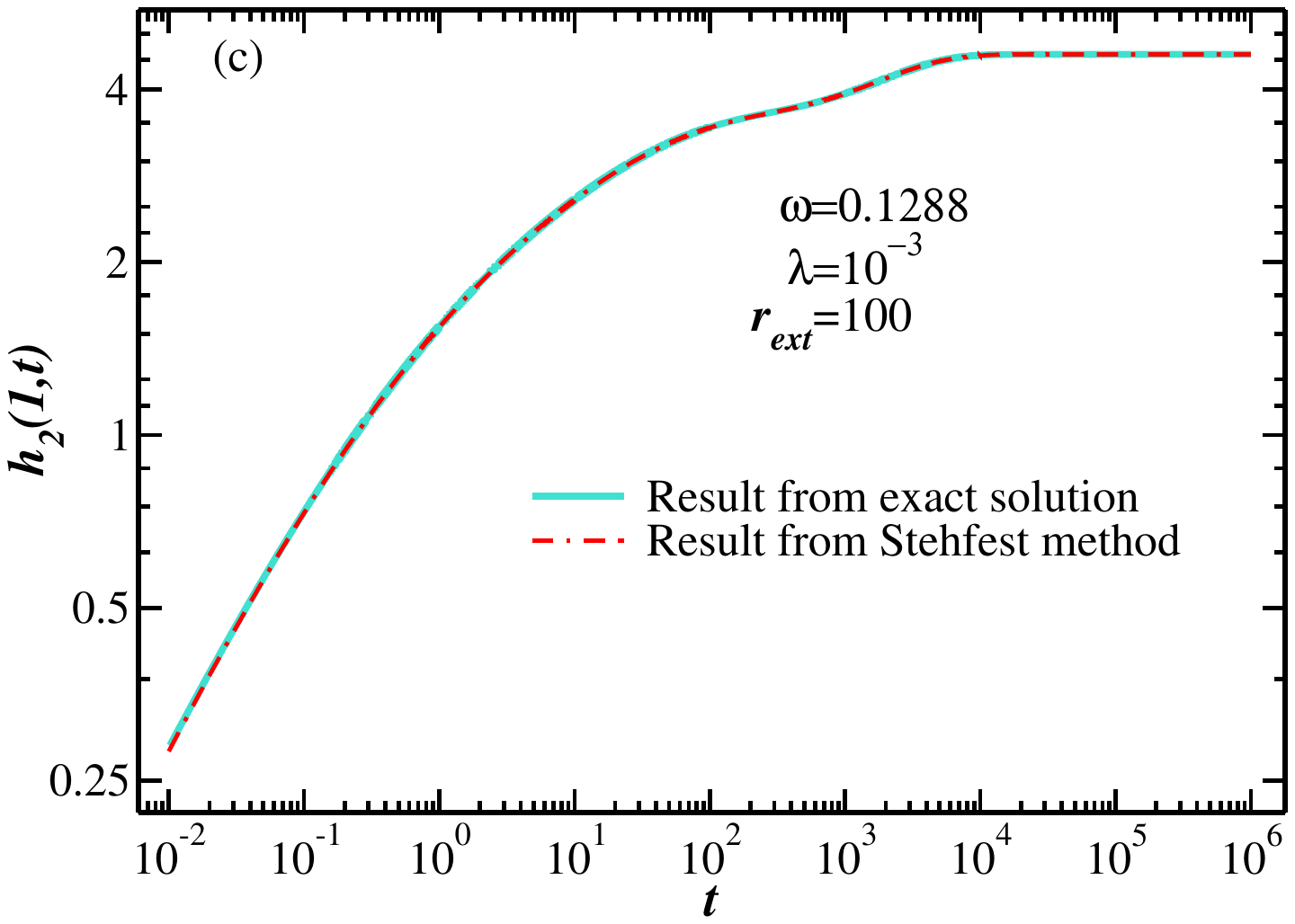}}
{\includegraphics[width=0.45\linewidth]{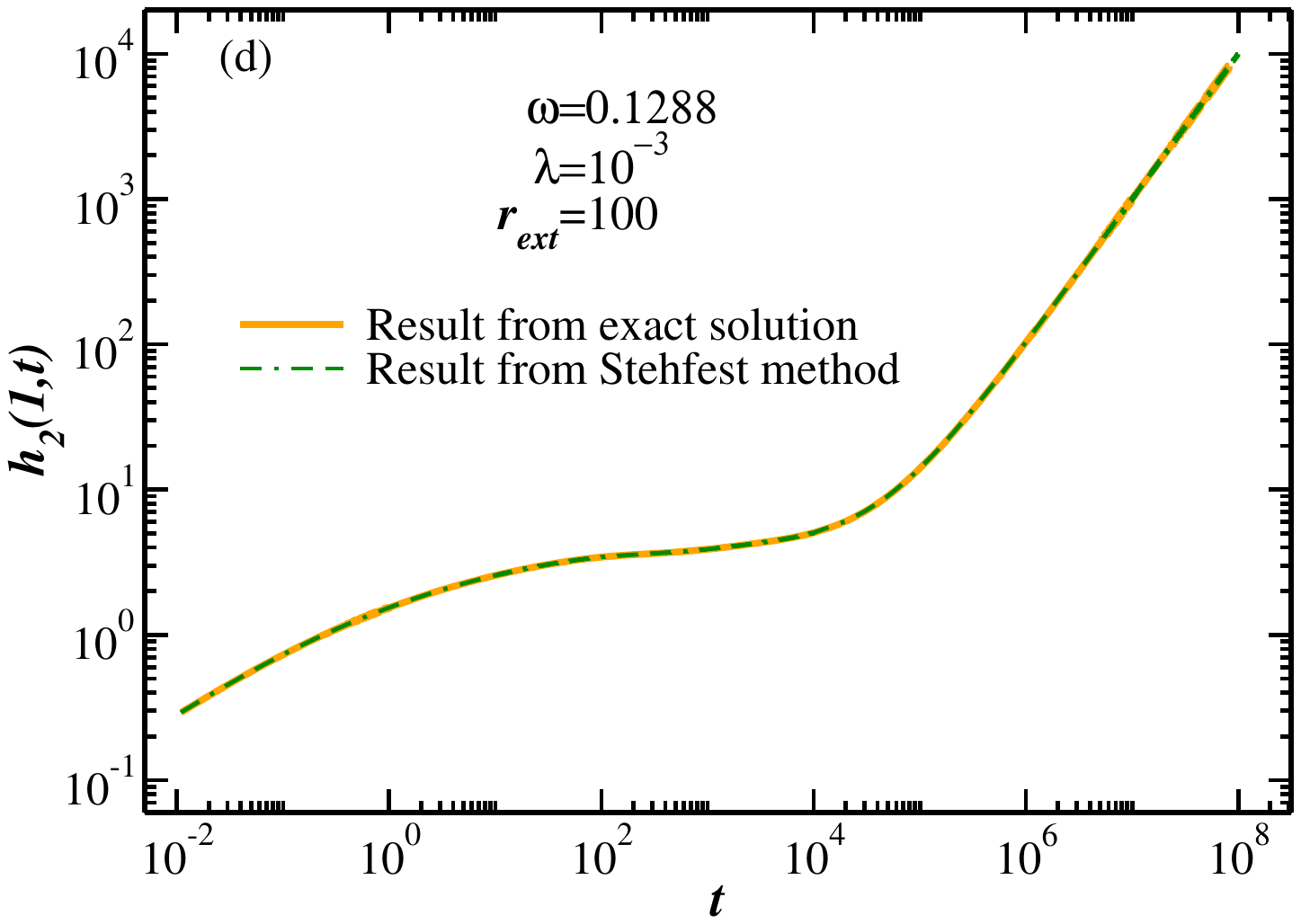}}\\
\caption{A comparison between the exact solutions from JLHT, and numerical solutions from the Stehfest method, is shown. It is observed a matching between both results. Eqs.~(\ref{solFlux_case1}), (\ref{Radialflux_case2}), (\ref{solWR_case3}), and (\ref{solWR_case4}), are used to obtain the results in the subfigures for the cases (a) DD-BCs, (b) DN-BCs, (c) ND-BCs, and (d) NN-BCs, respectively, while the results from inverse Laplace transform are obtained making use of Eqs.~(\ref{WR_j2D_caseA}), (\ref{WR_h2D_caseB-2}), (\ref{WR_h2D_caseC}), and (\ref{WR_h2D_caseD}). The parameters values used for the ``Ramp" rate function are $q_{\text{ext}} = 0.5$ and $\gamma = 10^{-3}$. In (a), (b), (c), and (d), are used $100$, $100$, $1000$, and $200$ $k_i$-roots, respectively 
}
\label{casoA_WR}
\end{figure}

Fig.~\ref{casoA_WR} shows the matching between the results from exact analytic solutions, Eqs.~(\ref{solFlux_case1}), (\ref{Radialflux_case2}), (\ref{solWR_case3}), and (\ref{solWR_case4}), and the results from the inversion of Eqs.~(\ref{WR_j2D_caseA}), (\ref{WR_h2D_caseB-2}), (\ref{WR_h2D_caseC}), and (\ref{WR_h2D_caseD}), respectively. Values of the parameters are given inside the frames of the figure. Figs.~\ref{casoA_WR}(a) and \ref{casoA_WR}(b) have graphs of flux into the wellbore due to fractures, while Figs.~\ref{casoA_WR}(c) and \ref{casoA_WR}(d) have graphs of pressure at the bottomhole. They exhibit the different production stages, which for Figs.~\ref{casoA_WR}(a) and \ref{casoA_WR}(b) are related to fractures, transition fractures-matrix, and matrix, and they are presented at early, middle, and long times, respectively~(\cite{warren1963behavior}). On the other hand, Figs.~\ref{casoA_WR}(b) and \ref{casoA_WR}(d), in addition, involve a stage dominated by recharge boundary effects, which are presented for a long-time production~(\cite{doublet1995decline,del2014pressure,wang2017transient}). Therefore, these latter figures include graphs with a transition matrix-recharge. In any case, the transitions between media with different permeabilities are around the inflection points~(\cite{uldrich1979method}). In Fig.~\ref{casoA_WR}(b), we note that the stage do\-mi\-na\-ted by the transition fracture-matrix is dimmed because the influx recharge effects arise close to this transition stage. In summary, the expected physical behavior of a fluid in a double-porosity medium is observed, as well as a perfect matching, to the naked eye, between the numerical and theoretical results. This latter point corroborates the exact analytical solutions presented in this work. It is worth mentioning that similar characteristic curves, as the described above, are obtained for other parameters choices, but for clarity are not shown in Fig.~\ref{casoA_WR}.

\subsection{Convergence of the Cinelli and exact solutions}

\begin{figure}[ht]
{\includegraphics[width=0.41\linewidth]{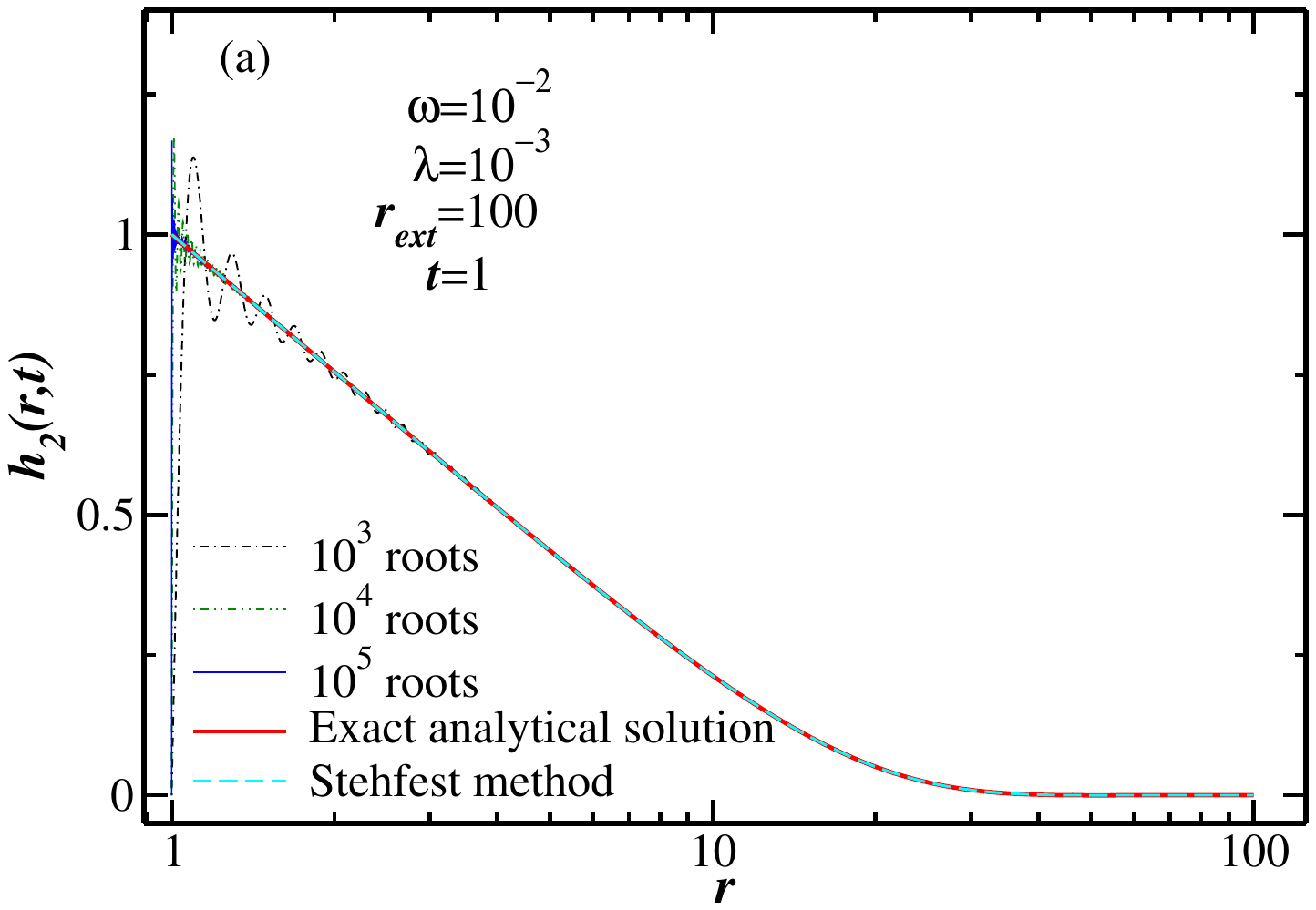}}
\put(-107.5,45){\includegraphics[width=2.8cm,height=2.8cm]{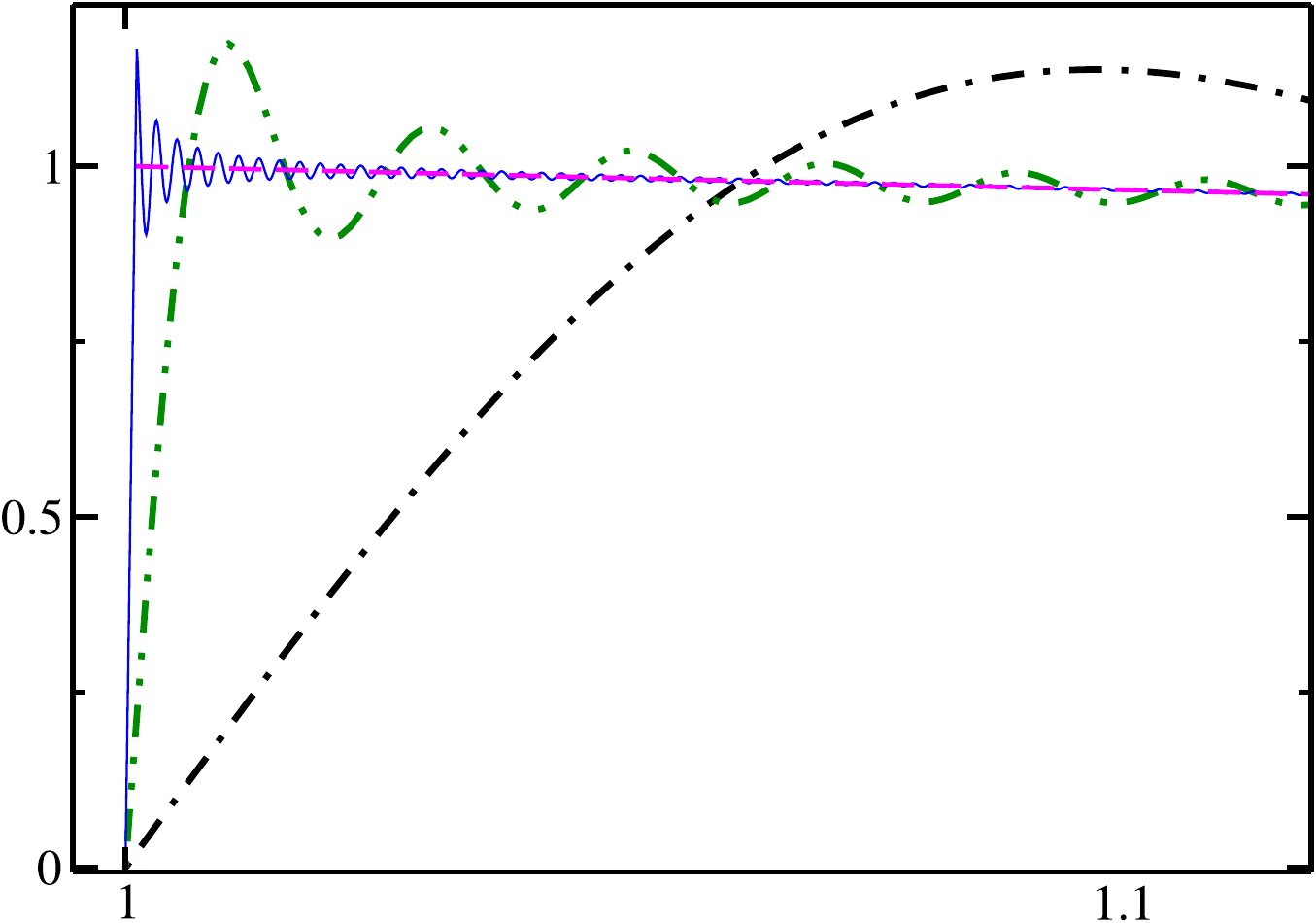}}
{\includegraphics[width=0.41\linewidth]{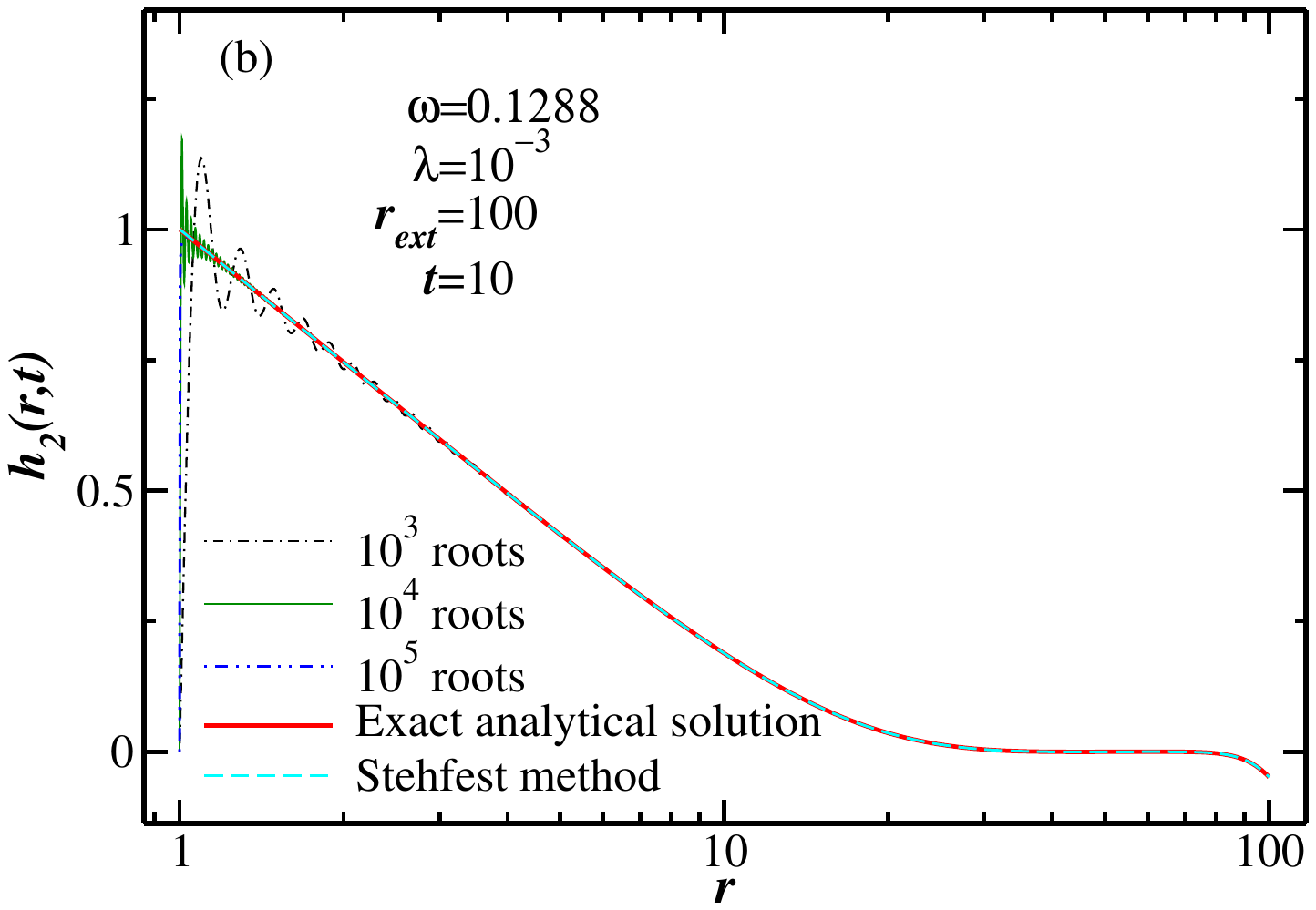}}
\put(-106.5,46){\includegraphics[width=2.8cm,height=2.8cm]{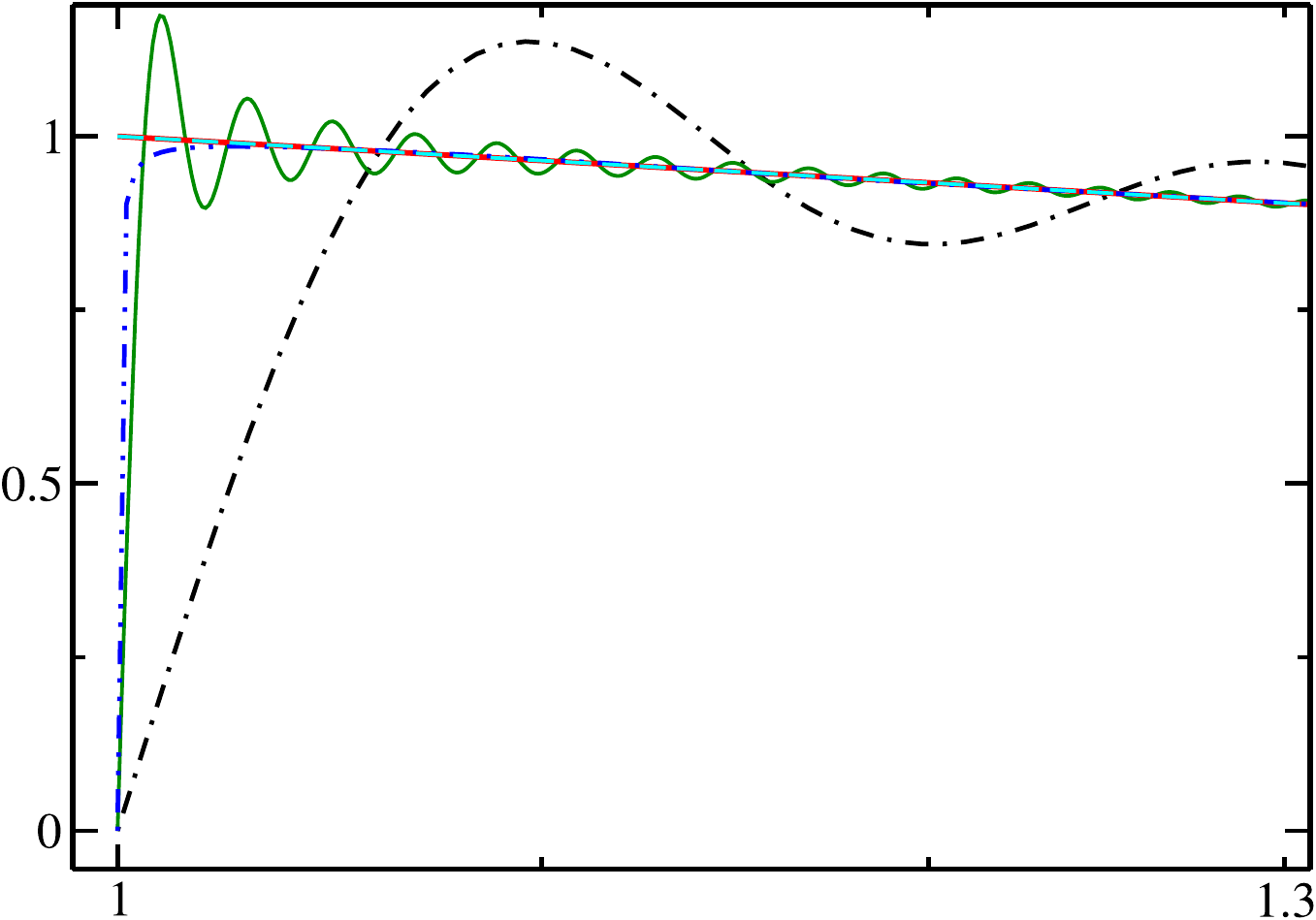}}\\
{\includegraphics[width=0.41\linewidth]{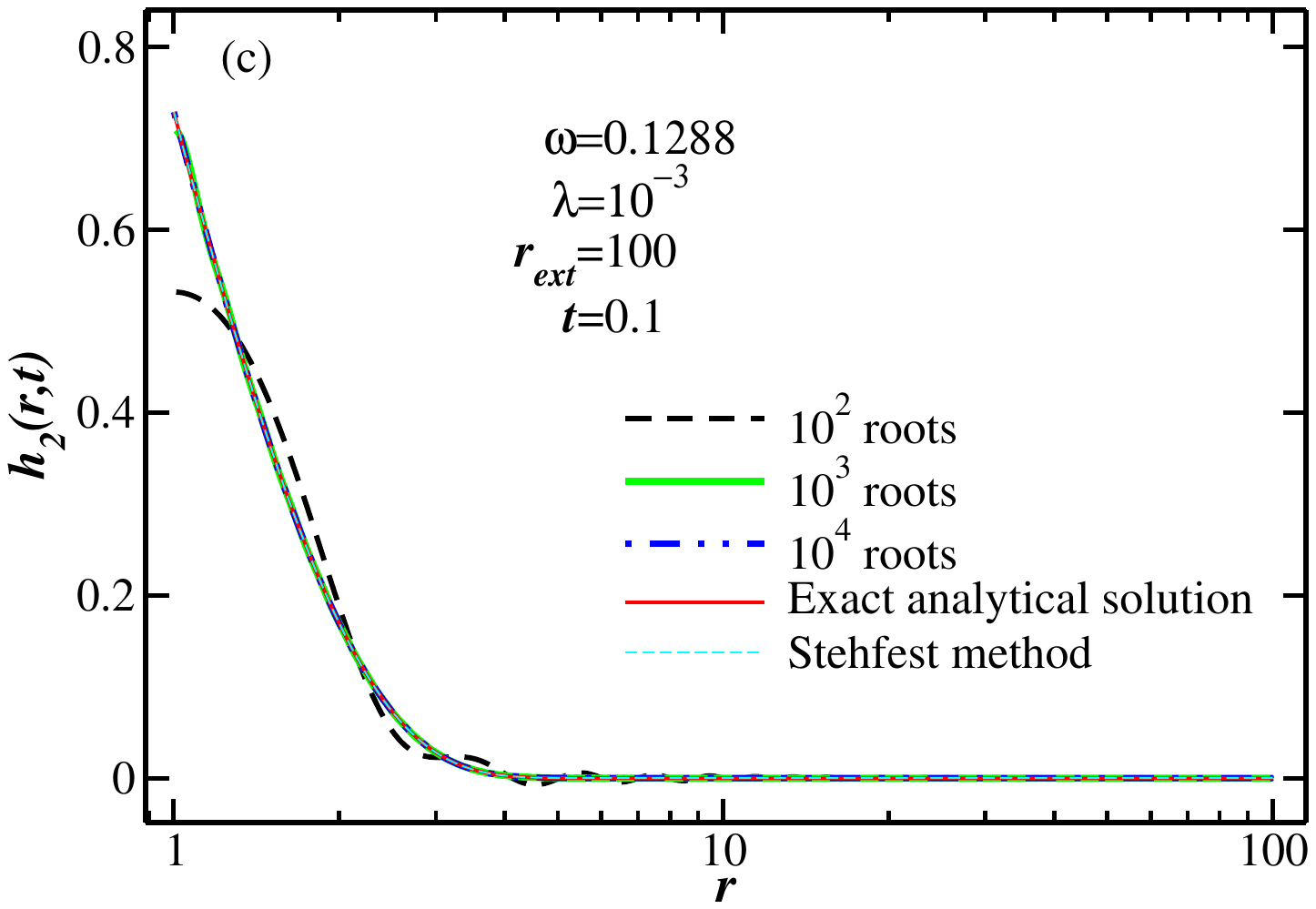}}
{\includegraphics[width=0.41\linewidth]{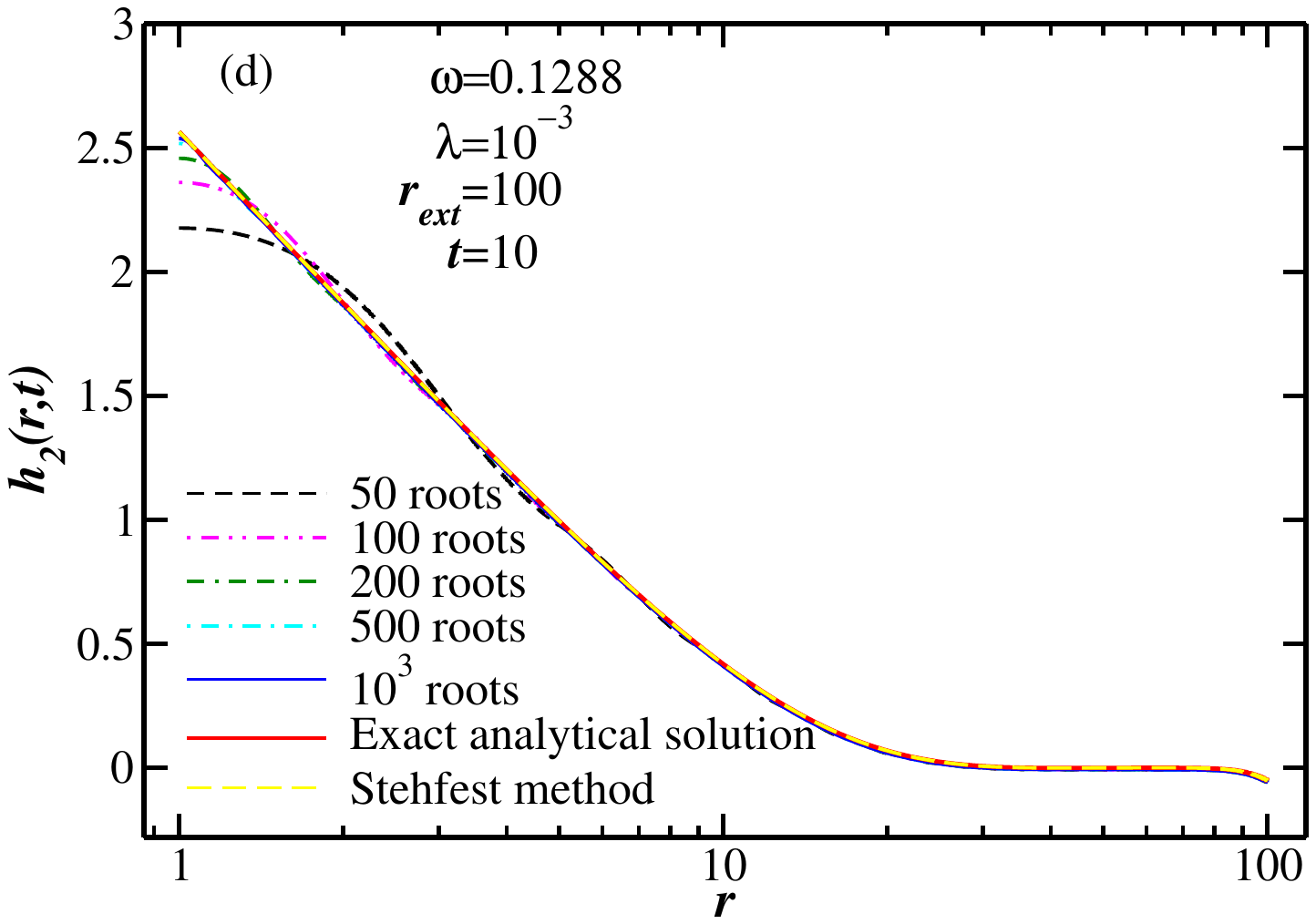}}
\put(-106.5,46){\includegraphics[width=2.8cm,height=2.8cm]{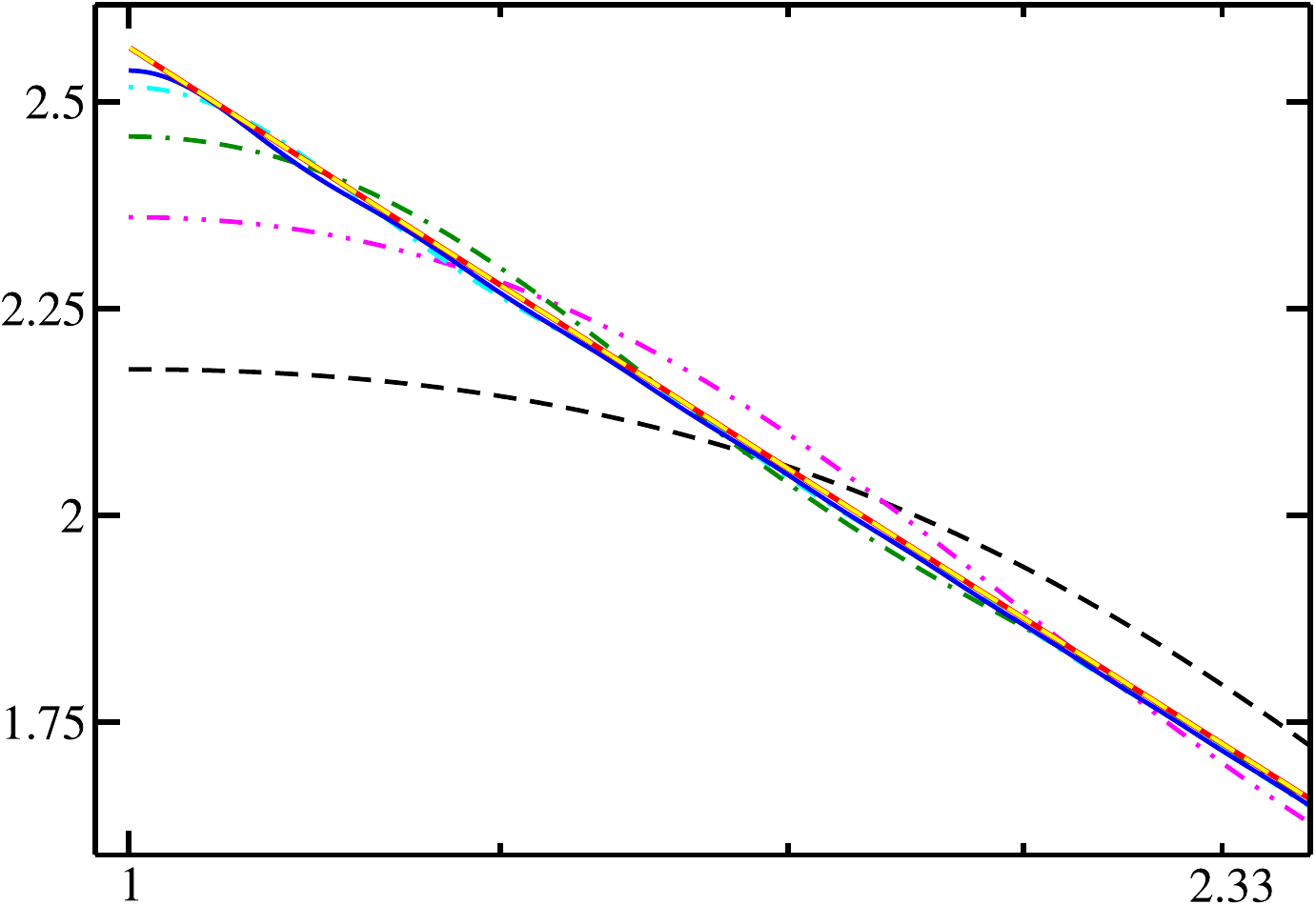}}\\
{\includegraphics[width=0.41\linewidth]{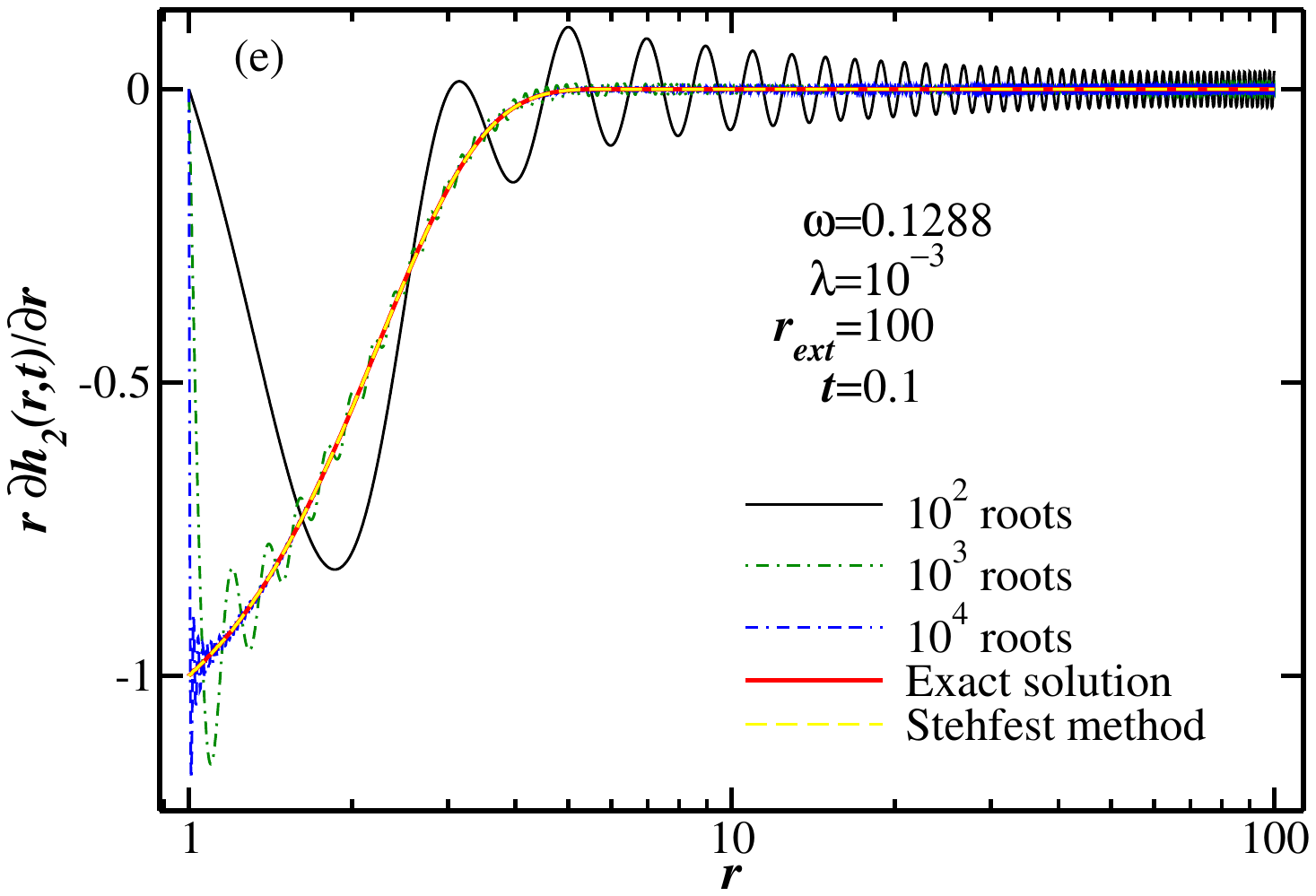}}
{\includegraphics[width=0.41\linewidth]{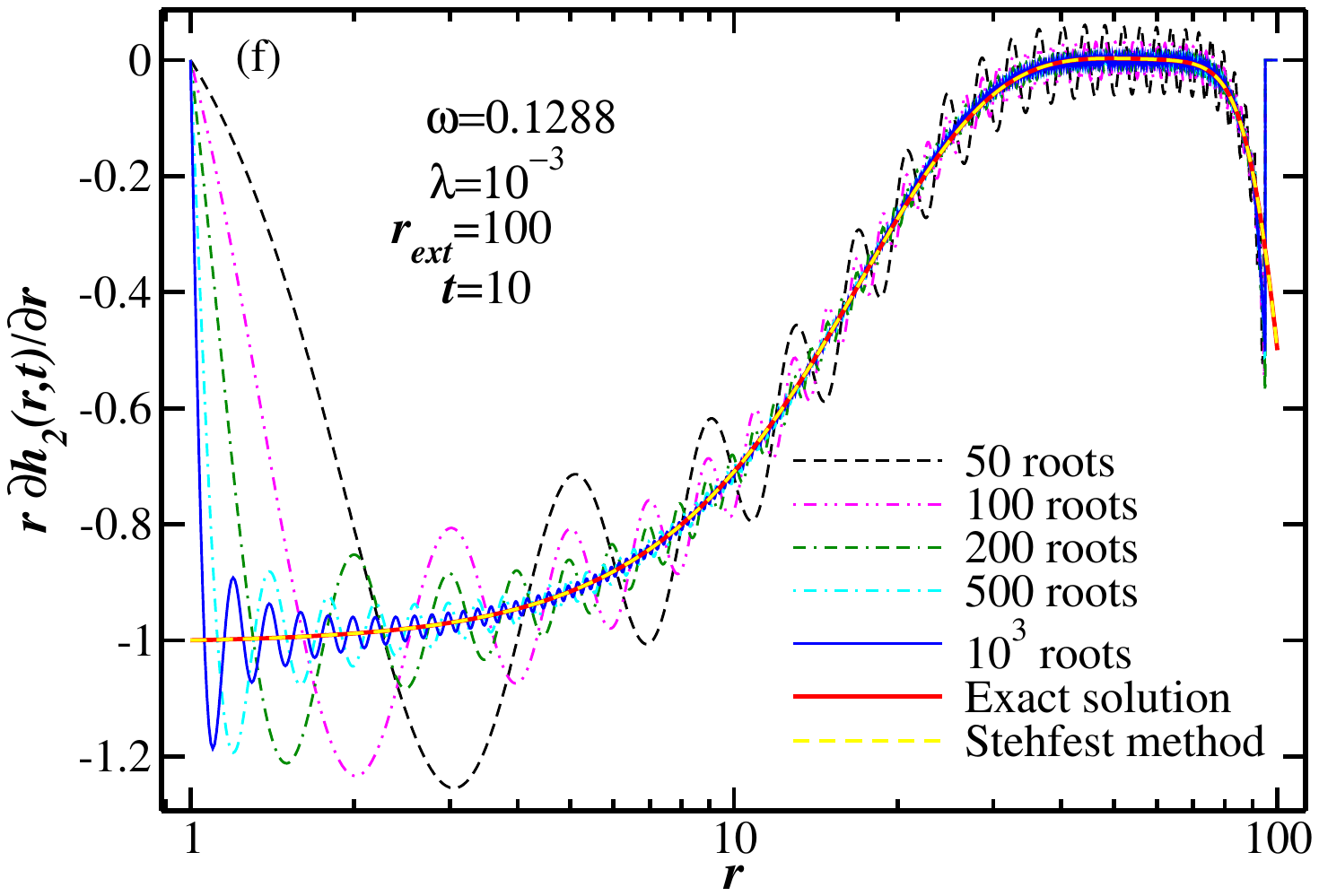}}\\
\caption{Results of the exact pressure (our $h_2(r,t)$ formulas) for the cases (a) DD-BCs, (b) DN-BCs, (c) ND-BCs, and (d) NN-BCs, and exact flux $j_2(r,t)$ for the cases (e) ND-BCs, and (f) NN-BCs. These results are compared with the corresponding ones obtained with Eqs.~(\ref{inverse_h_2(k,s)}) and with the Stehfest method applied to the equations of Table~\ref{TableAppendix}. It is remarkable how the results from Eqs.~(\ref{inverse_h_2(k,s)}), except for subfigures (c) and (d), show a slow converge and oscillations around the exact results, and also that the inner boundary condition is not hold in any case. For (c) and (d), the oscillatory behavior is diminished and the convergence is faster than in the other subfigures; note that in these cases the inner boundaries have flux conditions which never are meet. In all cases is observed a perfect matching between the exact solutions and the Stehfest method results
}
\label{case_oscillations}
\end{figure}

Graphs in Fig.~\ref{case_oscillations} contain the results obtained from applying Eqs.~(\ref{inverse_h_2(k,s)}) in our study model, i.e.~no closed formulas for the time independent infinite series are used. In this figure it is remarkable the systematic convergence of the solutions (towards the numerical result) by increasing the number of terms in Eqs.~(\ref{inverse_h_2(k,s)}). Also, notice that these equations hold a homogeneous condition in the inner boundary and that their convergence is very slow close to the bottomhole. In addition, note that the oscillations are increased by increasing the number of terms, as can be seen in Figs.~\ref{case_oscillations}~(a) and \ref{case_oscillations}~(b). These figures contain graphs of the DD-BCs and DN-BCs cases when $10^3$, $10^4$ and $10^5$ terms or roots are considered in the computations. Despite the great increase in the number of roots, the solutions oscillate considerably around the numerical inversion. Note that these latter numerical results match with the ones from Eqs.~(\ref{solWR_case1}) and (\ref{solWR_case2}) and are free of oscillations. On the other hand, Figs.~\ref{case_oscillations}~(c) and \ref{case_oscillations}~(d) contain graphs for ND-BCs and NN-BCs cases. Again, we observed a systematic convergence by increasing the number of terms. In fact, for ND-BCs and NN-BCs cases, the oscillation quickly diminishes when more terms of the series are considered in computation. Even so, the number of terms is large compared to the number of terms used in the computation by means of equations that involves the closed formulas, Eqs.~(\ref{solWR_case3}) and (\ref{solWR_case4}). For example, in Fig.~\ref{casoA_WR} hundreds of terms are used, while in Fig.~\ref{case_oscillations} thousands. However, the behavior like the one given for Figs.~\ref{case_oscillations}~(a) and \ref{case_oscillations}~(b) is found when flux along reservoir (at fixed time) is plotted, see Figs.~\ref{case_oscillations}~(e) and \ref{case_oscillations}~(f). This is because these graphs show the convergence of solution toward the value of the inner BC. It is worth remarking that results in Fig.~\ref{case_oscillations} are computed either with the solutions of Cinelli or the exact ones, and that both have infinite sums. Therefore, the absence of oscillations in the exact solutions is due to the fact that the infinite sums (corresponding to the stationary limit) could be approximated with simplified analytical expressions, i.e.~the oscillations are attributable to computations with a finite number of terms of the time-independent series. In summary, it is remarkable that: 1) direct use of Eqs.~(\ref{inverse_h_2(k,s)}) leads to the correct results, except at inner boundary, and 2) when the closed formulas are used, the inhomogeneous BCs hold independently of the number of terms.

\subsection{Existence of stationary solutions}

As we have highlighted, the solution of
\begin{linenomath*}
\begin{eqnarray}\label{eqODELaplace}
\frac{1}{r}\frac{\partial}{\partial r}\left( r \frac{\partial h_{2,s}}{\partial r} \right)= 0, \;\;\;\;\mbox{which is } h_{2,s}(r) = A \log (r) + B,
\end{eqnarray}
\end{linenomath*}
is involved in Eq.~(\ref{eqReplaceClosed}) and this is stationary. For this reason, the solution of model (\ref{2Pmodel}) is represented as the sum of a (time-independent) stationary solution and a (time-dependent) transitory solution. This latter solution has the same form as Cinelli formulas, replacing $\tilde{h}_{2}(k_i,t)$ in Eqs.~(\ref{inverse_h_2(k,s)}) by $\tilde{g}(k_i,t)$. Therefore, the transitory part holds by default for homogeneous BC and, for this reason, the stationary part $h_{2,s}(r)$ should describe the inner BC.

\begin{figure}[ht]
{\includegraphics[width=0.5\linewidth]{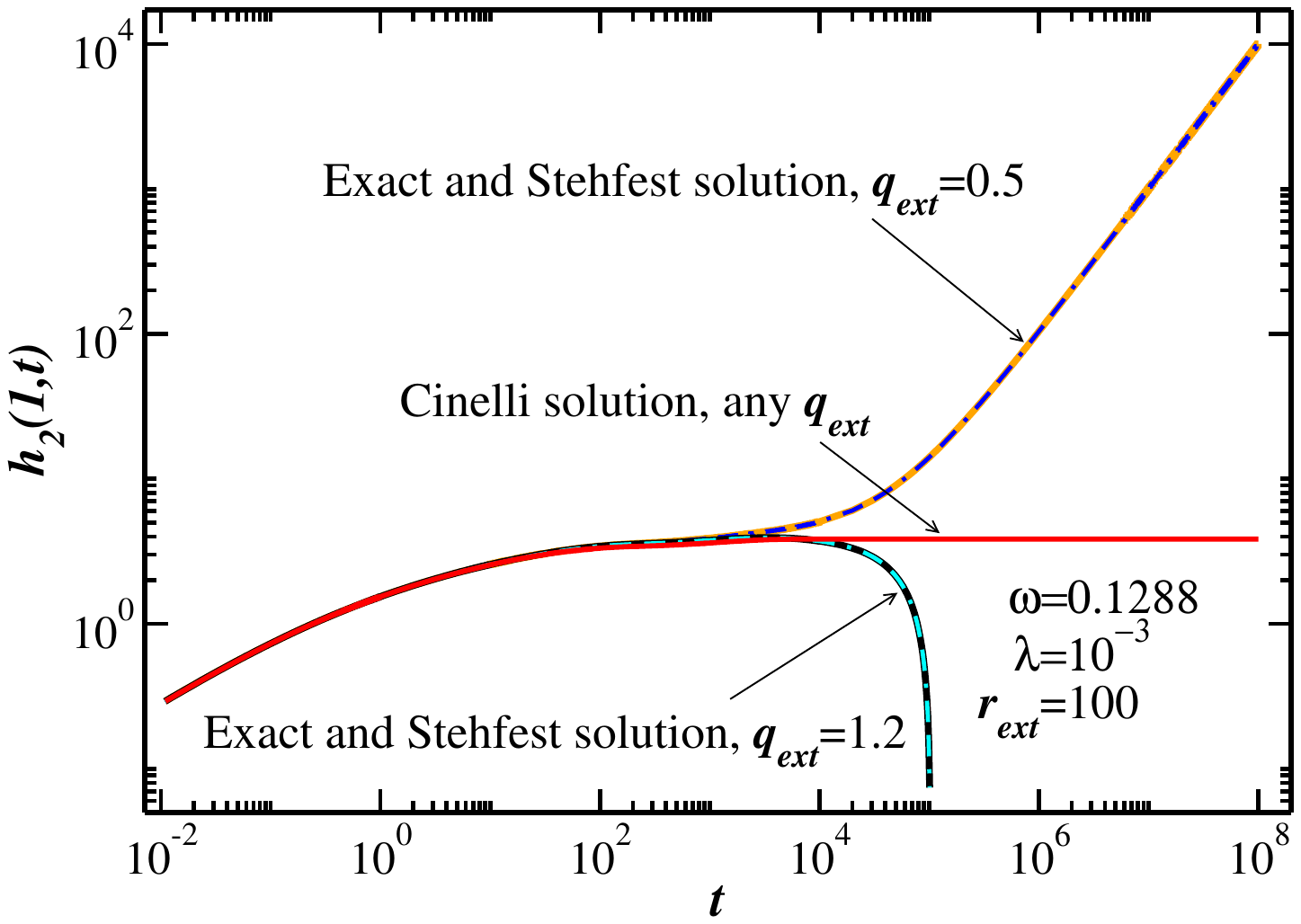}}
\caption{Drawdown pressure curves for the NN-BCs case. The flattened line is the result from the Cinelli formula or Eq.~(\ref{inverse_h_2(k,s)d}), and the dashed lines are obtained using the exact result given by Eq.~(\ref{solWR_case4}). For comparison, the inverted numerical results from Eq.~(\ref{WR_h2D_caseD}) are also shown in dashed lines. In all cases $3\times 10^4$ terms are used. The considered $\gamma$ value is $10^{-3}$
}
\label{case_Cinelli}
\end{figure}

Because the non-uniqueness condition required for a well-posedness boundary value problem, the Eq.~(\ref{eqODELaplace}) has no stationary solution when a fixed flow is imposed at the inner boundary with an influx through the outer boundary. Solution methods can lead to uniqueness by specifying a conservation principle or splitting the problem into well-posedness problems, see, for example, exercises 18.3.10 (c) and 18.3.19 in~\citet{greenberg1988advanced}.  Regarding this remark, stationary solutions arise when the inner and outer boundary conditions are equal, i.e.~when $A$ has the same value in both boundaries. Thus, the flux has a long-time stationary behavior because no net flow enters or leaves the reservoir. Under these conditions, we have validated that our methodology, in relation with Eq.~(\ref{eqReplaceClosed}), works well. Otherwise, when the total flux increases or decreases at any time, a dynamic behavior predominates over time. In this case, we have observed that the formula in~\citet{cinelli1965extension} for NN-BCs has inconsistencies when the influx recharge has a predominant effect in the system. In fact, the flattened line in Fig.~\ref{case_Cinelli} exhibits such mistake, since the results given by this formula remain constant for long-time, for any value of $q_{\text{ext}}$ chosen. This behavior disagrees with the one obtained using the Stehfest method, in which there is an increase in the asymptotic behaviors when $q_{\text{ext}}<1$ and a decrease when $q_{\text{ext}}>1$. These characteristics are consistent with those observed in Refs.~(\cite{van1949application,del2014pressure}) and, as we can see in Fig.~\ref{case_Cinelli}, they are reproduced by our exact solution, Eq.~(\ref{solWR_case4}), which is valid for any value of $q_{\text{ext}}$ considered.

It is worth mentioning that we have added deliberately the time-dependent terms in the right-hand side of Eq.~(\ref{solWR_case4}) to obtain the correct behavior; these terms are identified by being outside of the infinite series indicated there. As it was stated in Section~\ref{Section_Hankel-Laplace_Sol}, these added terms are obtained by expanding Eq.~(\ref{WR_h2D_caseD}) about $s=0$, and, subsequently,  taking the inverse Laplace transform to the result of this expansion. It can be seen, that these terms are not considered in the Cinelli relationships, Eq.~(\ref{inverse_h_2(k,s)d}), and their inclusion leads to a correct matching with the numerical results, as can be seen in Fig.~\ref{case_Cinelli}. Similarly, it can be proved that stationary solutions for DD-BCs, DN-BCs and ND-BCs cases can be recovered using this procedure.

\section{Characteristic behavior of reservoirs with influx recharge}\label{Section_characteristic_behavior_of_solutions}

In this section, characteristic curves of drawdown pressure and flux for reservoirs with influx recharge are presented. Drawdown pressure curves are analyzed in the framework of the pressure derivative of Bourdet~(\cite{bourdet1983new,bourdet1989use}). Subsequently, a comparative analysis between models with NN-BCs is done, considering a model of a single-porosity and non-zero influx and a second model of double-porosity but with a closed-boundary. This analysis is done with the purpose of elucidating whether they can be quantitatively equivalent.

It is worth mentioning, before beginning the analysis, that for naturally fractured reservoirs, usual values of of the interporosity flow coefficient are $10^{-10} \leq \lambda \leq 10^{-4}$, while the fracture storage coefficient $\omega$ is in  the order of $10^{-3}$ or $10^{-2}$~(\cite{bourdet2002well}). Regarding the work of~\citet{kuhlman2015multiporosity}, we have $\omega \lesssim$ 0.1\%. The parameter $q_{\text{ext}}$ (or $\gamma$) has not been characterized because the scarcity of studies in this direction. On the other hand, according to the Warren and Root model, the parameter values admit the following ranges: $0< \omega \leq1$, $\lambda>0$, $q_{\text{ext}} \geq 0$, and $\gamma > 0$. In a similar way to other works~(\cite{mavor1979transient,wang2017transient}), we use these latter ranges to show the characteristic behavior of the study model. Indeed, best-fitting curves of the Warren and Root model could lead to a large value of $\omega$, e.g.~see~\citet{camacho2014well}.

\subsection{Characteristic curves}\label{Section_behaviors}

\begin{figure}[ht]
{\includegraphics[width=0.45\linewidth]{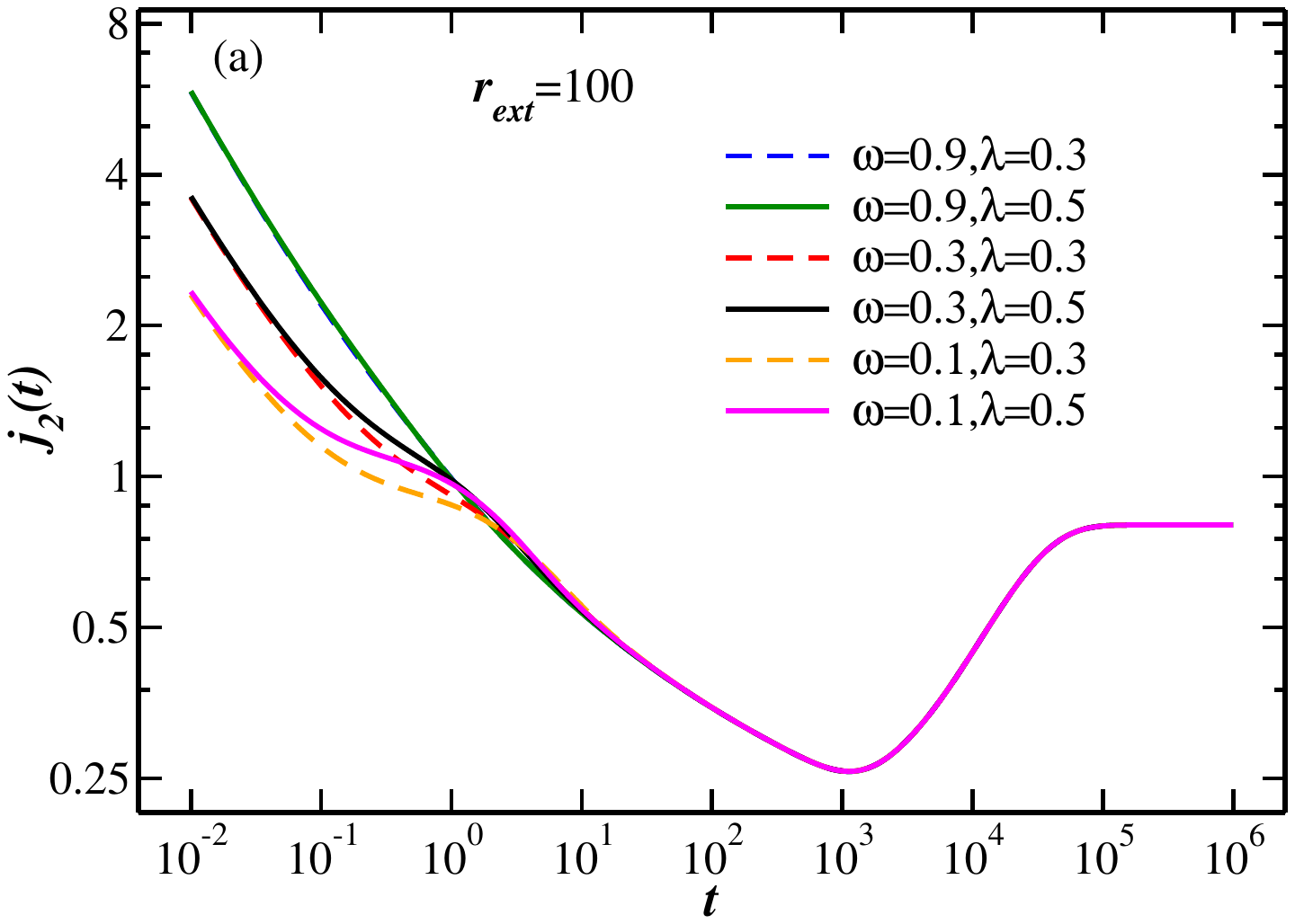}}
{\includegraphics[width=0.45\linewidth]{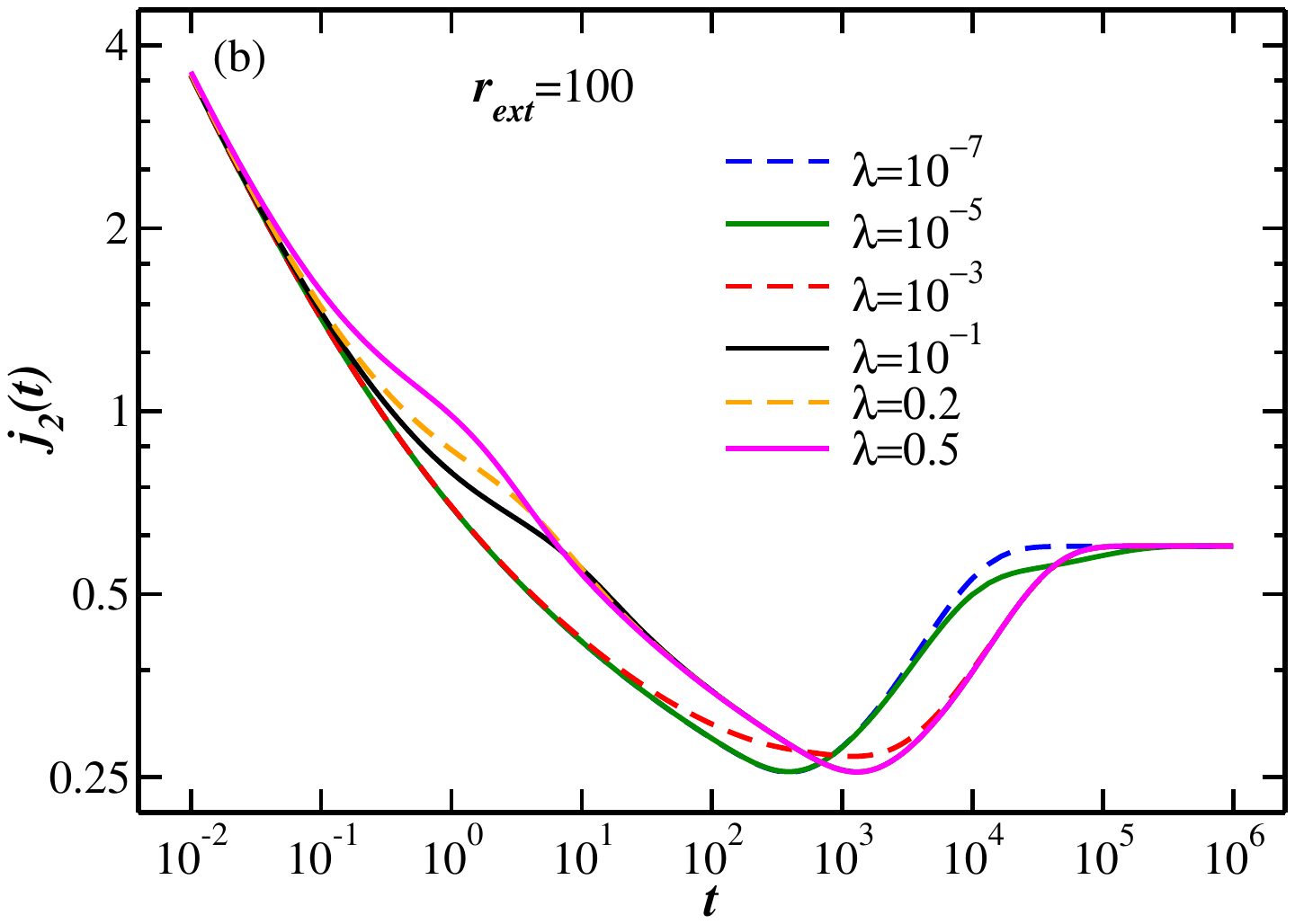}}\\
{\includegraphics[width=0.45\linewidth]{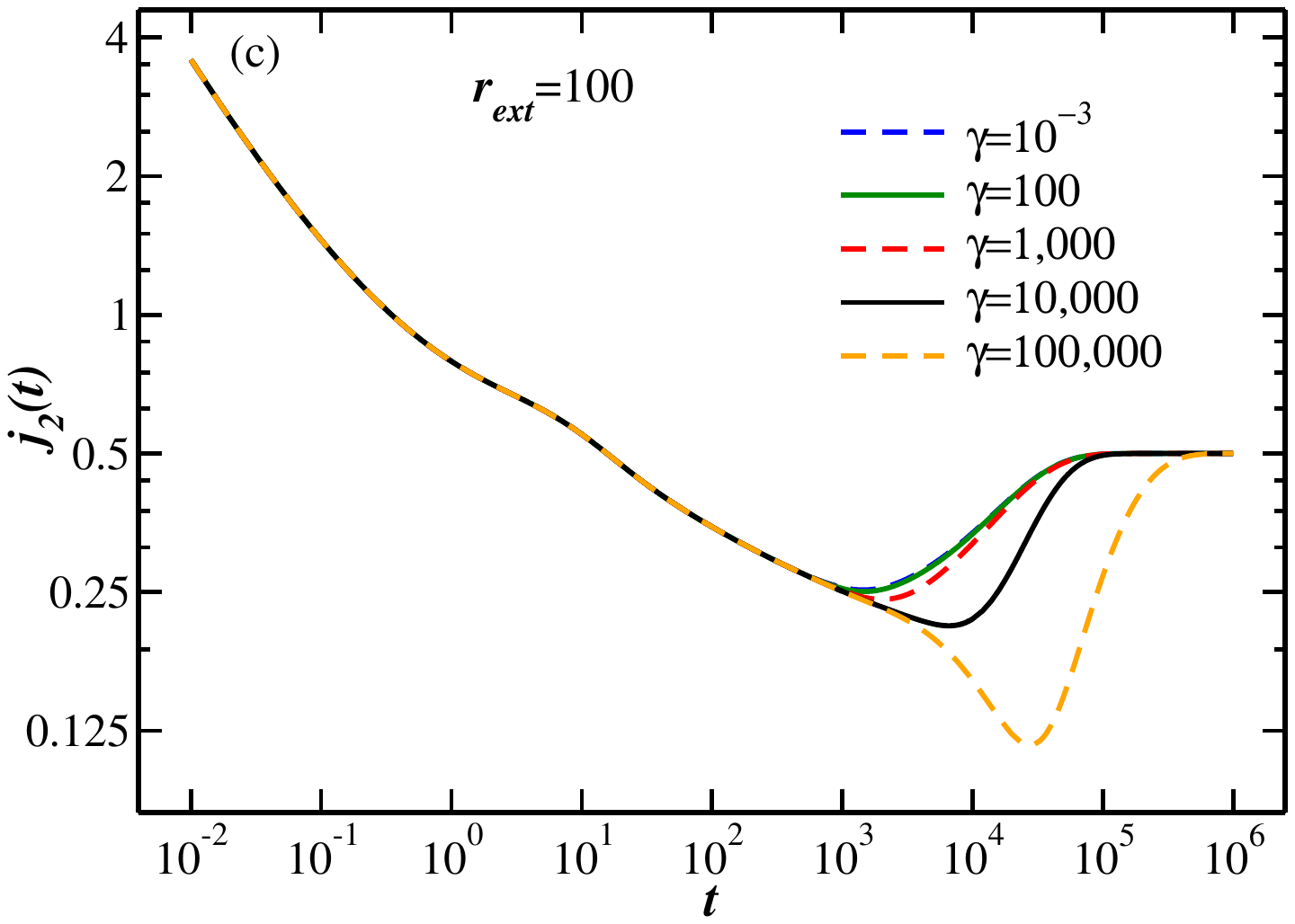}}
{\includegraphics[width=0.45\linewidth]{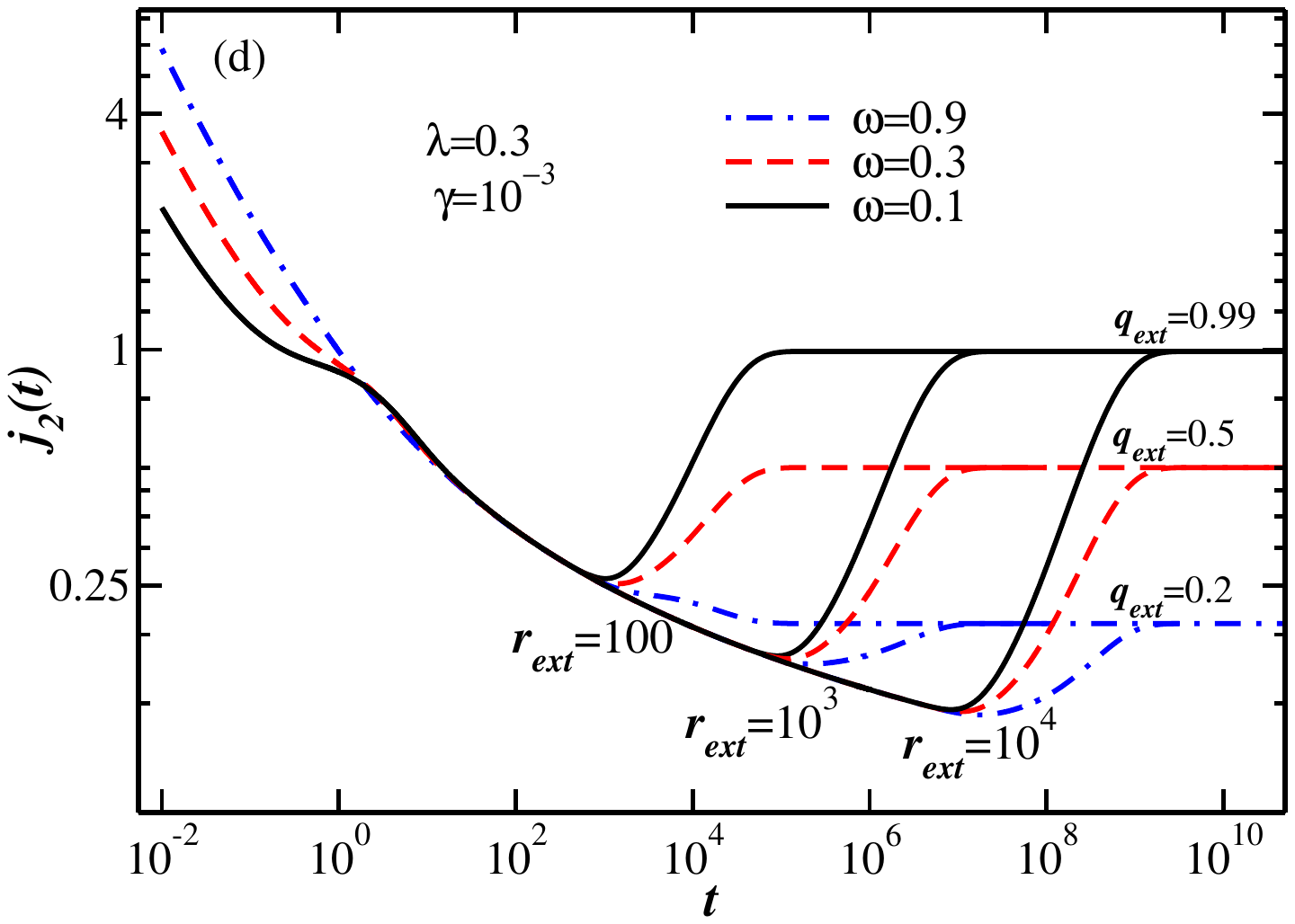}}
\caption{Flux characteristic curves obtained from Eq.~(\ref{Radialflux_case2}) for the DN-BCs case. The frames show the effect of the following parameters on the curves: (a) $\omega$ and $\lambda$, (b) $\lambda$, (c) $\gamma$, and (d) $\omega$, $r_{\text{ext}}$, and $q_{\text{ext}}$. The rest of the parameters values, in addition to the ones given in frames, are: (a) $q_{\text{ext}}=0.8$ and $\gamma=10^{-3}$, (b) $\omega=0.3$, $q_{\text{ext}}=0.6$ and $\gamma=10^{-3}$, (c) $\omega=0.3$, $q_{\text{ext}}=0.5$ and $\lambda=0.1$ 
}
\label{characteristic_behavior_DN-BC}
\end{figure}

\subsubsection{DN-BCs case}\label{Section_DN-BCs_case}

Taking into account that the condition imposed at the bottomhole is constant pressure, Eq.~(\ref{Radialflux_case2}) leads to the flux behaviors shown in Fig.~\ref{characteristic_behavior_DN-BC}. Therein, the differences between the graphs come from varying the values of the storage $\omega$, the interporosity flux coefficient $\lambda$, the slope of the ``Ramp" $\gamma$, the outer influx factor $q_{\text{ext}}$, and the outer radius $r_{\text{ext}}$. Note that the interporosity flow coefficient is a  first-order mass transfer coefficient~(\cite{goltz1991analytical,moench1995convergent,pedretti2014apparent}). In Fig.~\ref{characteristic_behavior_DN-BC} (a) is evident that at early times $\omega$ has a major effect on the flux. Namely, a small $\omega$ leads to a transition fracture-matrix more pronounced, while the graph for $\omega=0.9$ does not show the characteristic form of such transition. Fig.~\ref{characteristic_behavior_DN-BC} (a) also exhibits that, for fixed $\omega$, this form is affected by changing the $\lambda$ value.  The $\lambda$ effect in the flux can also be seen in Fig.~\ref{characteristic_behavior_DN-BC} (b); for $\lambda\geq 10^{-1}$ transitions with a negative half-slope are observed, while the graphs for the rest of values of $\lambda$ are similar to the ones from a single-porosity medium. On the other hand, since $\gamma$ is an influx parameter, its effect is seen in Fig.~\ref{characteristic_behavior_DN-BC} (c) for a long time production. Because a large $\gamma$ implies a slow influx recharge into the reservoir, the flux declination is more pronounced when the $\gamma$ value is increased. In addition, the effect of $r_{\text{ext}}$ becomes clear in Fig.~\ref{characteristic_behavior_DN-BC} (d). Namely, the characteristic curves show that a larger $r_{\text{ext}}$ leads to a greater flux drop. In this figure, we give three graphs for each value of $r_{\text{ext}}$, each of these depending on a couple of values $q_{\text{ext}}$ and $\omega$ as indicated there. Notice that the flux drop is recovered at long time when the flux becomes stationary with a $q_{\text{ext}}$ value; it can be seen for every graph in Fig.~\ref{characteristic_behavior_DN-BC}. We do not include solutions for a closed reservoir, $q_{\text{ext}}=0$, however, for this case the flux tends to zero and the solutions have not minimums as the showed in Fig.~\ref{characteristic_behavior_DN-BC}.

\subsubsection{NN-BCs case}\label{Section_NN-BCs_case}

\begin{figure}[th]
{\includegraphics[width=0.45\linewidth]{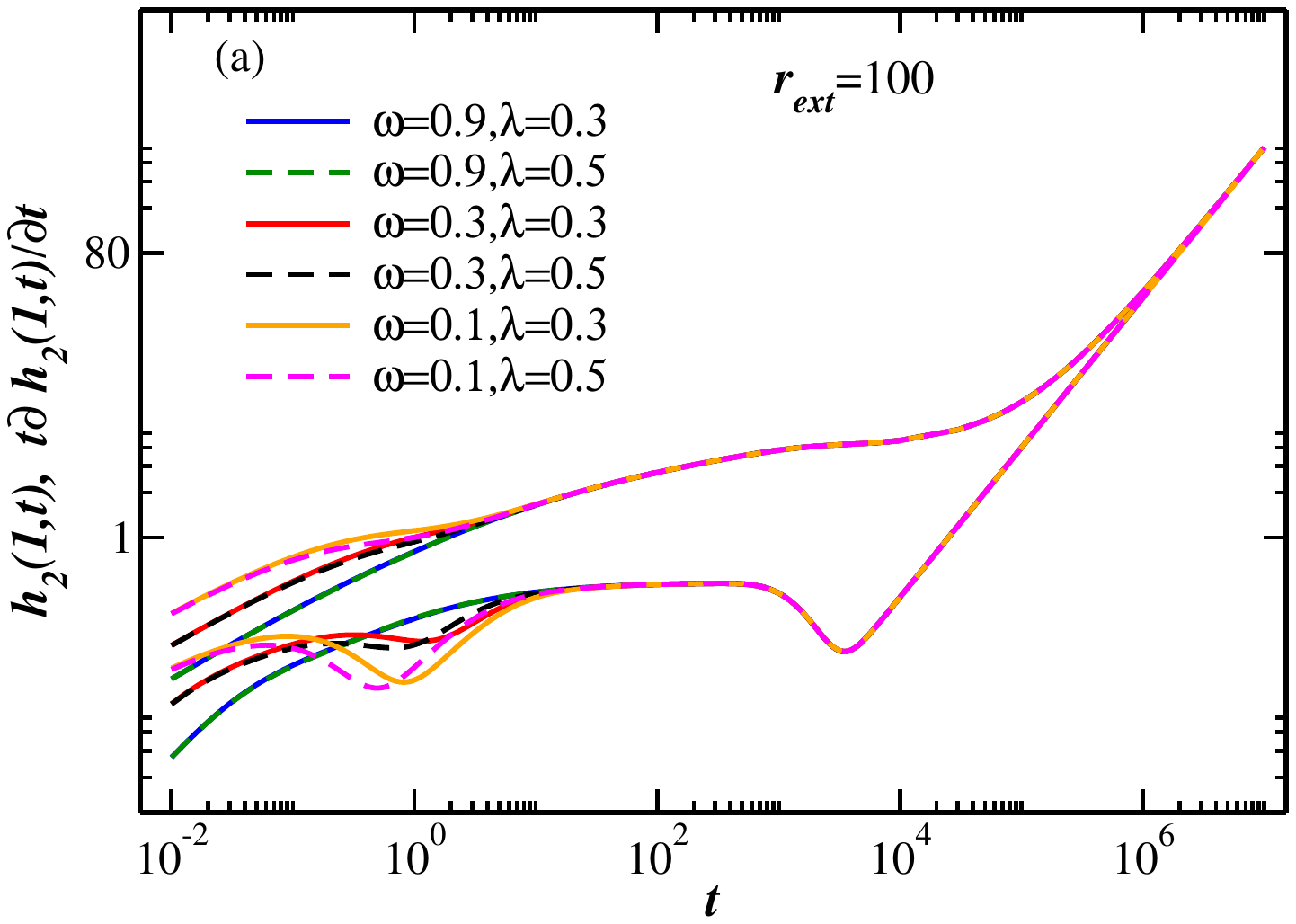}}
{\includegraphics[width=0.45\linewidth]{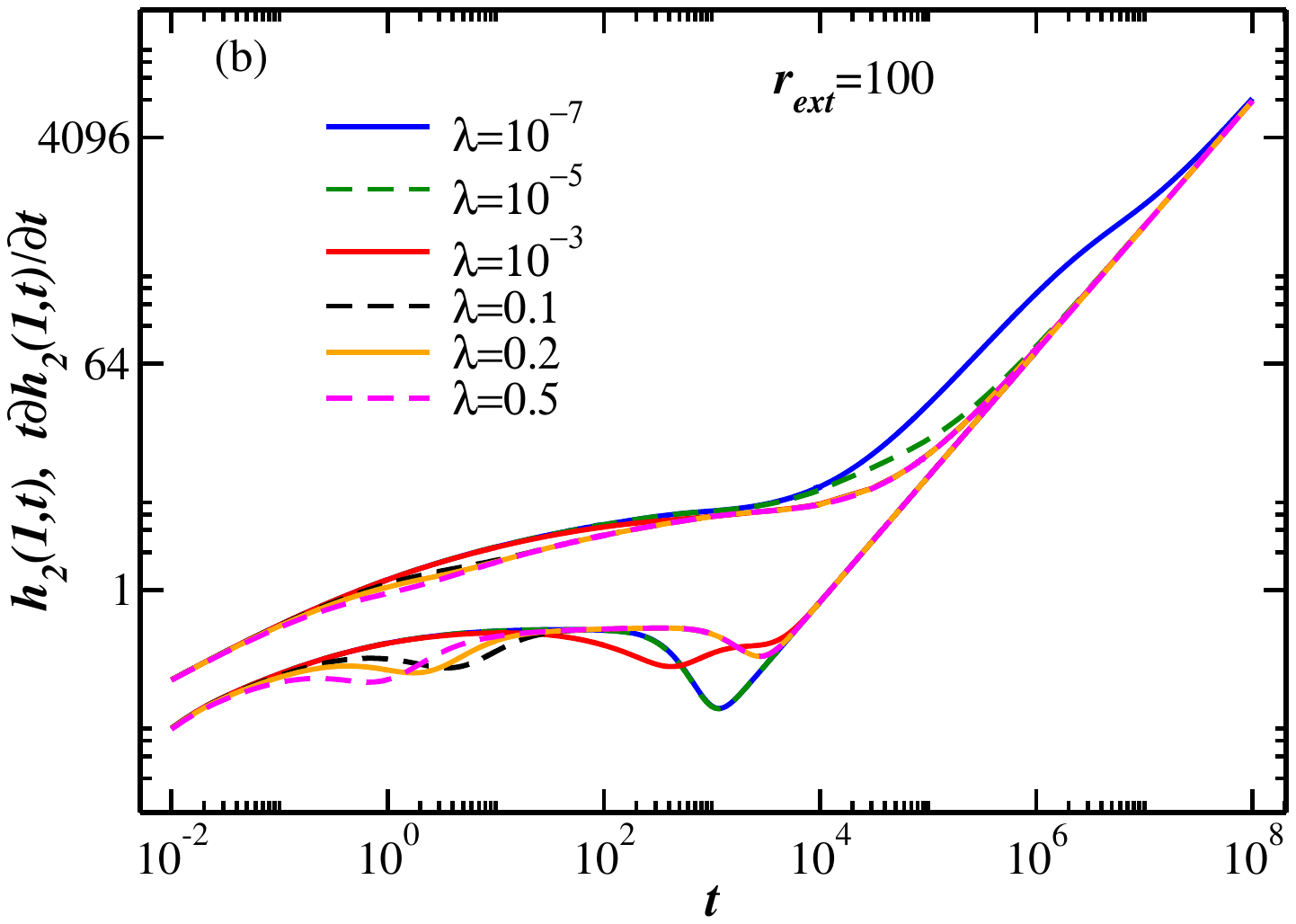}}\\
{\includegraphics[width=0.45\linewidth]{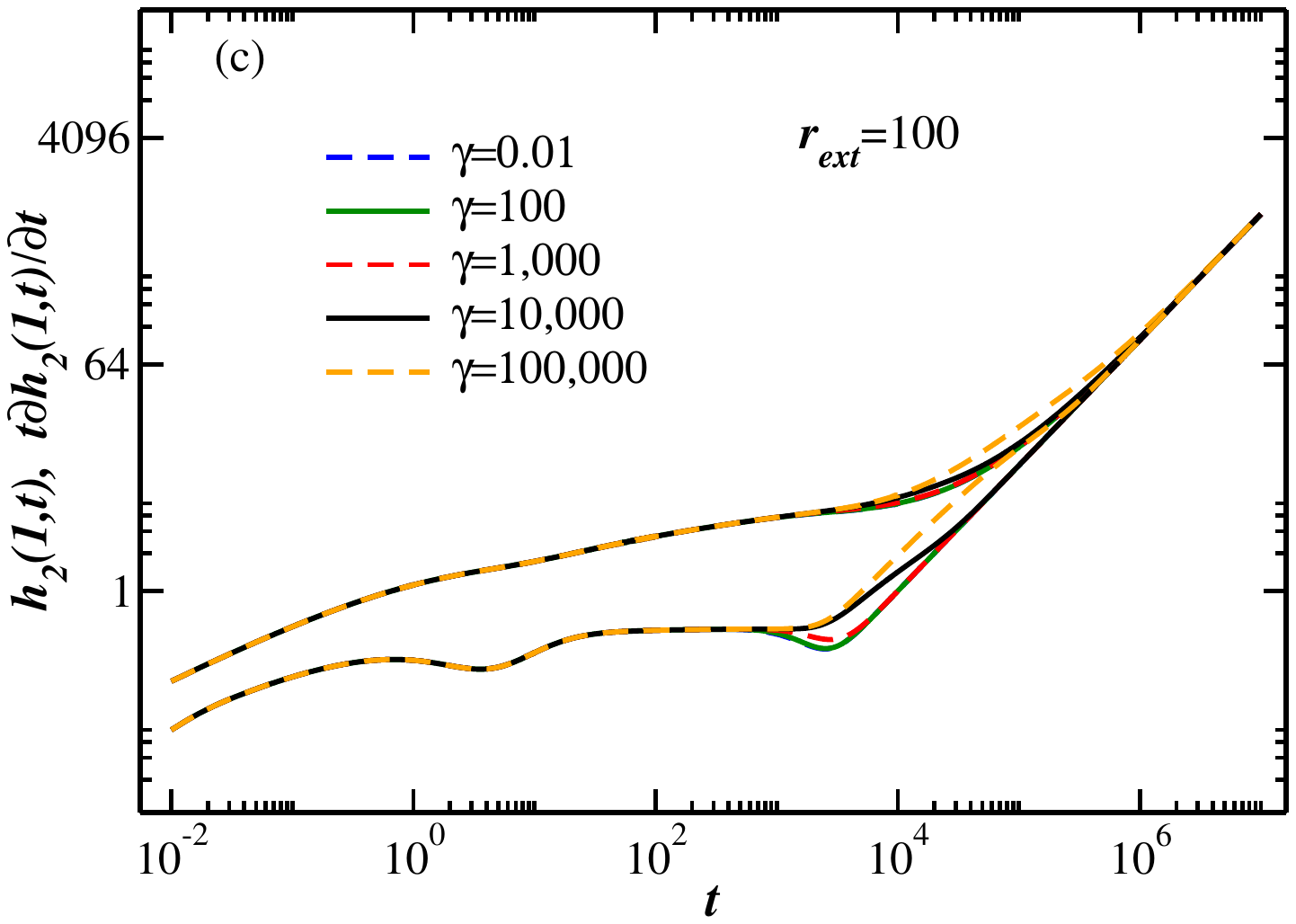}}
{\includegraphics[width=0.45\linewidth]{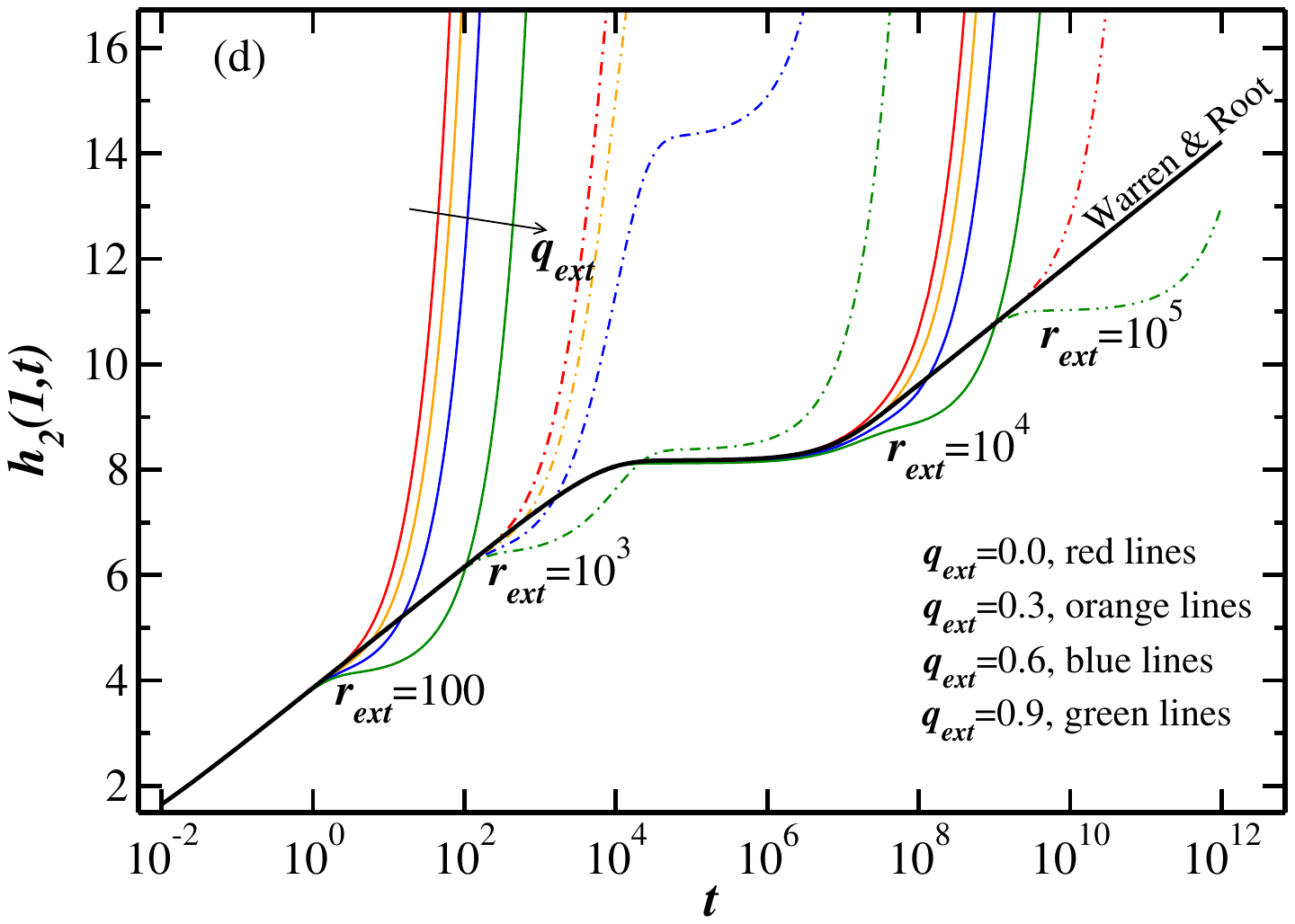}}
\caption{Drawdown pressure characteristic curves, Eq.~(\ref{solWR_case4}), and their Bourdet derivatives for the NN-BCs case. The frames show the effect of the following parameters on the curves: (a) $\omega$ and $\lambda$, (b) $\lambda$, (c) $\gamma$, and (d) $r_{\text{ext}}$ and $q_{\text{ext}}$.  The rest of the parameter values, in addition to the ones given in frames, are: (a) $q_{\text{ext}}=0.8$ and $\gamma=10^{-3}$, (b) $\omega=0.3$, $q_{\text{ext}}=0.6$ and $\gamma=10^{-3}$, (c) $\omega=0.3$, $q_{\text{ext}}=0.5$ and $\lambda=0.1$, and (d) $\omega=0.001$, $\lambda=10^{-7}$, and $\gamma=10^{-3}$ 
}
\label{characteristic_behavior_NN-BC}
\end{figure}

Examples of drawdown pressures and their Bourdet derivatives are found in Fig.~\ref{characteristic_behavior_NN-BC}. We can see how changes in $\omega$ have effect for early time, the effect by varying $\lambda$ covers the solution domain, while changes due to $\gamma$ are presented for long time; see Figs.~\ref{characteristic_behavior_NN-BC} (a), (b), and (c), respectively. Similar characteristic behaviors have been observed in other works~(\cite{del2014pressure,wang2017transient}), however, these works were not for double-porosity media. The drawdown pressure show several stages as the time evolve, which are influenced by fractures, transition fractures-matrix, matrix, transition matrix-influx, and influx. Furthermore, the Bourdet derivatives have minimums in correspondence with the transition stages. It can be noticed that both characteristic curves have the same linear behavior at long times, with a slope value of 1. However, note that this fact leads to an extension of the range of the matrix stage when $\lambda$ is very small, as can be seen for $\lambda = 10^{-7}$ in Fig.~\ref{characteristic_behavior_NN-BC} (b). A similar situation occurs in Fig.~\ref{characteristic_behavior_NN-BC} (c), where $\gamma=10^5$ implies a longer time for the recharge to take effect. In addition, the effect of $r_{\text{ext}}$ becomes clear in Fig.~\ref{characteristic_behavior_NN-BC} (d). Namely, the characteristic curves show that a high $r_{\text{ext}}$ and $t$ the drawdown behaviors are that of the Warren \& Root model. In this figure, we give a set of graphs for each value of $r_{\text{ext}}$ and different values of $q_{\text{ext}}$, each of these sets is indicated by different types of line. Notice also that the drawdown behavior is that of a closed reservoir, except when $r_{\text{ext}} \rightarrow \infty$, i..e.~the behavior at long time of a closed reservoir, or $q_{\text{ext}}=0$, is similar to that of a reservoir with influx recharge, when $0 < q_{\text{ext}} < 1$.

As mentioned above, the drawdown pressure and its Bourdet derivative have the same asymptotic curve at long times. Indeed, this is obtained from Eq.~(\ref{solWR_case4}) and its derivative,
\begin{linenomath*}
\begin{equation}\label{eqlongtimes}
t\frac{\partial h_2}{\partial t} = \left( \frac{2}{r_{\text{ext}}^2}-\frac{2q_{\text{ext}}}{r_{\text{ext}}^2-1} \right)t = h_2.
\end{equation}
\end{linenomath*}
In the previous result, we use the fact that the transitory part of solution tends to zero when time increases, as we deduced before in the discussion of Fig.~\ref{case_Cinelli}. The influence of the parameter $q_{\text{ext}}$ on pressure drawdown curves is analyzed in Refs.~(\cite{del2014pressure}) and (\cite{wang2017transient}), where a single-porosity and triple-porosity reservoir are studied, respectively. In these works it is observed that the influx recharge has no influence on the early and middle production stages. Furthermore, it is remarkable the following points: 1) when $q_{\text{ext}} = 1$ implies $\partial h_2/\partial t = 0$, i.e.~there is a stationary solution with a constant pressure at the outer boundary; 2) when $q_{\text{ext}}<1$ there is a pseudo-steady-state solution and the influence of $q_{\text{ext}}$ is similar to the case when there is a closed boundary; and 3) when $q_{\text{ext}}>1$ the pressure of the reservoir increases, while the bottomhole pressure decreases with a negative slope.

\begin{figure}[ht]
{\includegraphics[width=0.45\linewidth]{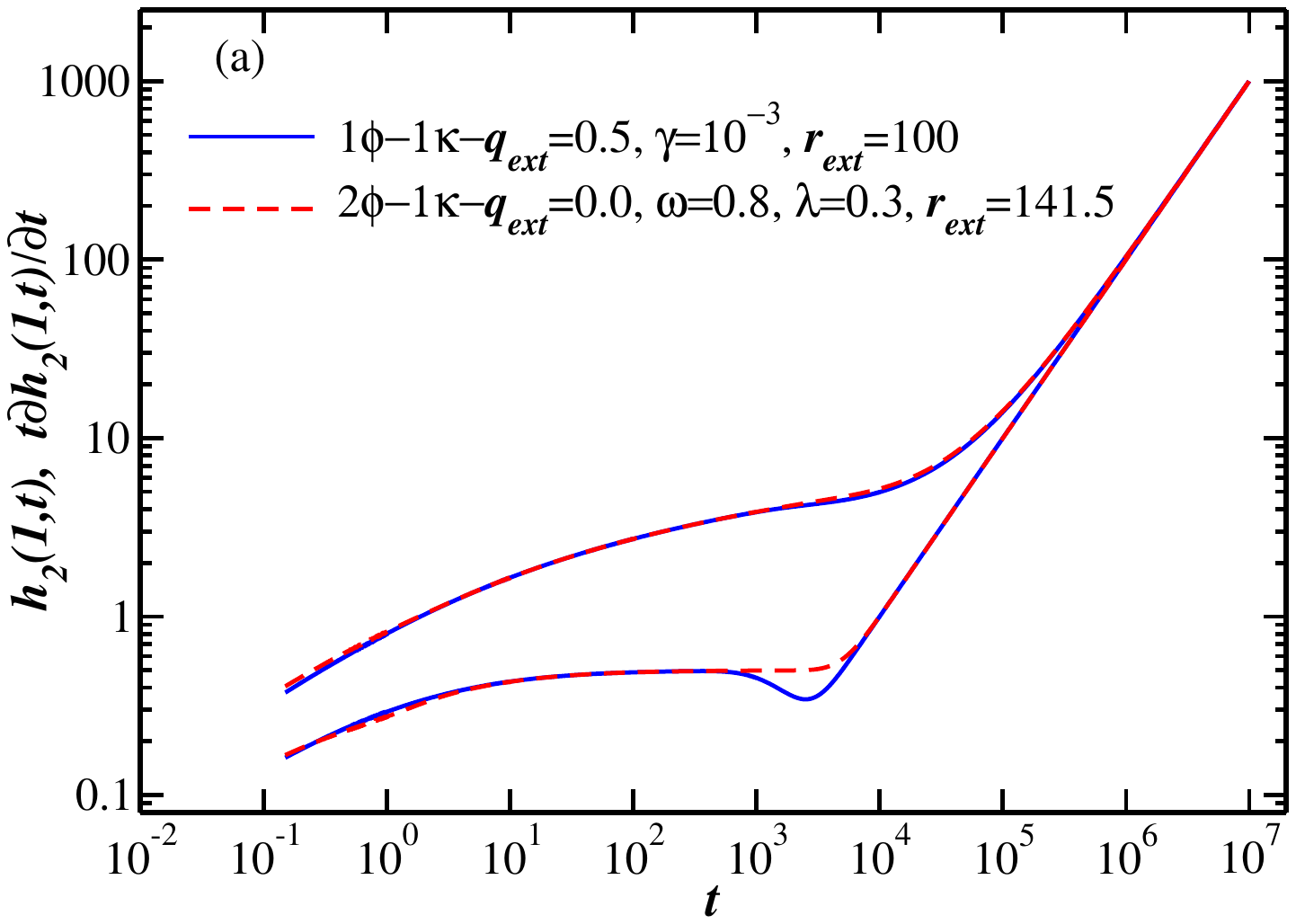}}
{\includegraphics[width=0.45\linewidth]{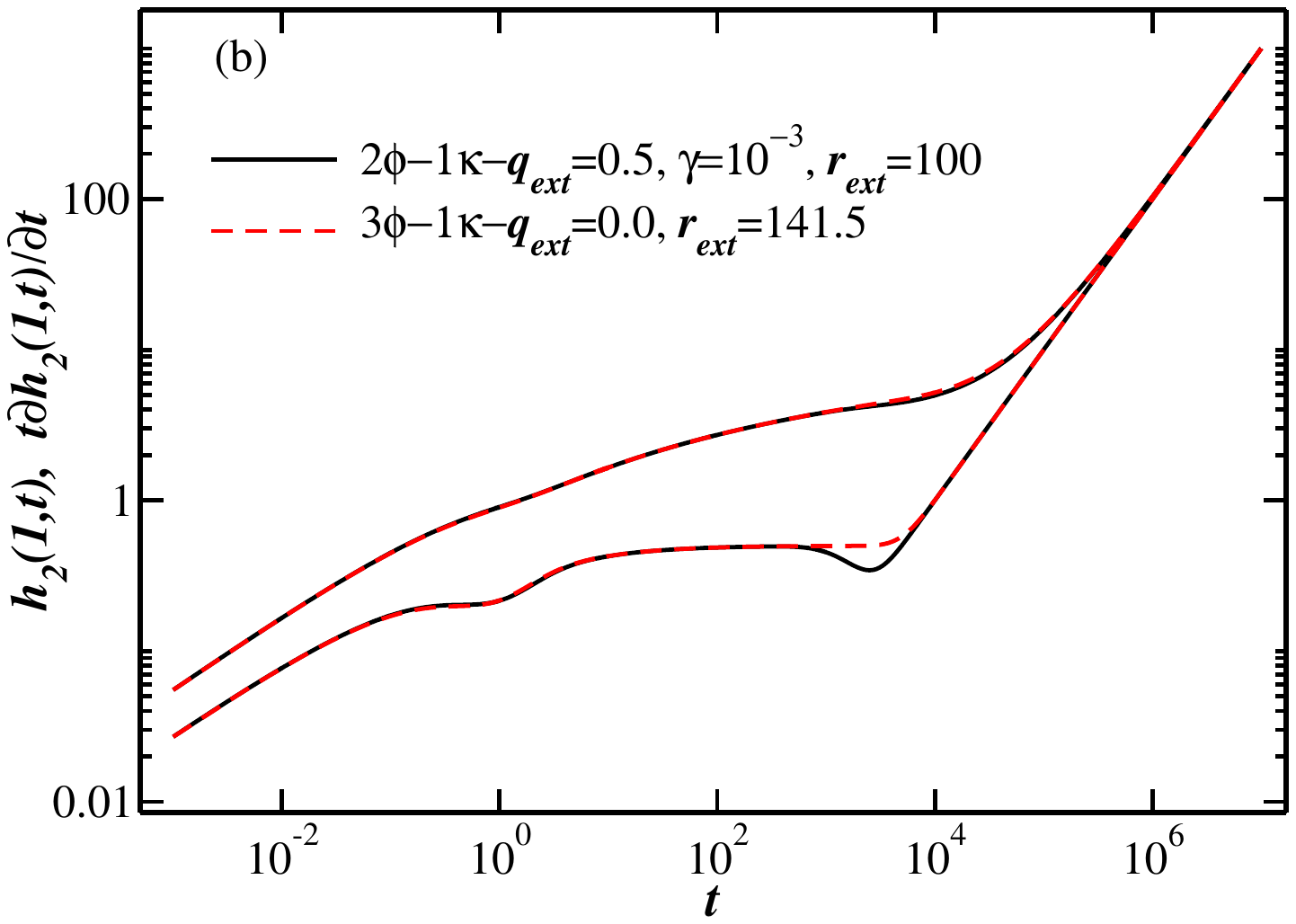}}\\
\caption{Comparison of the drawdown pressure and its Bourdet derivative of the following models: (a) $1\phi$-$1\kappa$-$q_{\text{ext}}\neq 0$ vs $2\phi$-$1\kappa$-$q_{\text{ext}}=0$ and (b) $2\phi$-$1\kappa$-$q_{\text{ext}}\neq 0$ vs $3\phi$-$1\kappa$-$q_{\text{ext}}=0$ (this triple-porosity model considers the differential equations (1)-(3) in~\citet{wu2007triple} together with the NN-BCs given in this work). As we can see in each frame, the difference between the drawdown curves are difficult to distinguish one of another. However, their Bourdet derivative show a clear difference in the transition period toward the outer boundary flow regime, i.e.~in this transition period there is a minimum that depends on $q_{\text{ext}}$. The graph of the triple-porosity model is obtained with $\omega_v=0.05$, $\omega_f=0.4$, $\lambda_{fv}=0.005$, $\lambda_{fm}=0.5474$, and $\lambda_{vm}=0.1643$ ($v$, $f$, and $m$ are for vuggs, fractures, and matrix, respectively), while the graph of the double porosity model with influx recharge is obtained with $\omega=0.4$ and $\lambda=0.5$
}
\label{curves_with_or_without_recharge}
\end{figure}

\begin{figure}[ht]
{\includegraphics[width=0.45\linewidth]{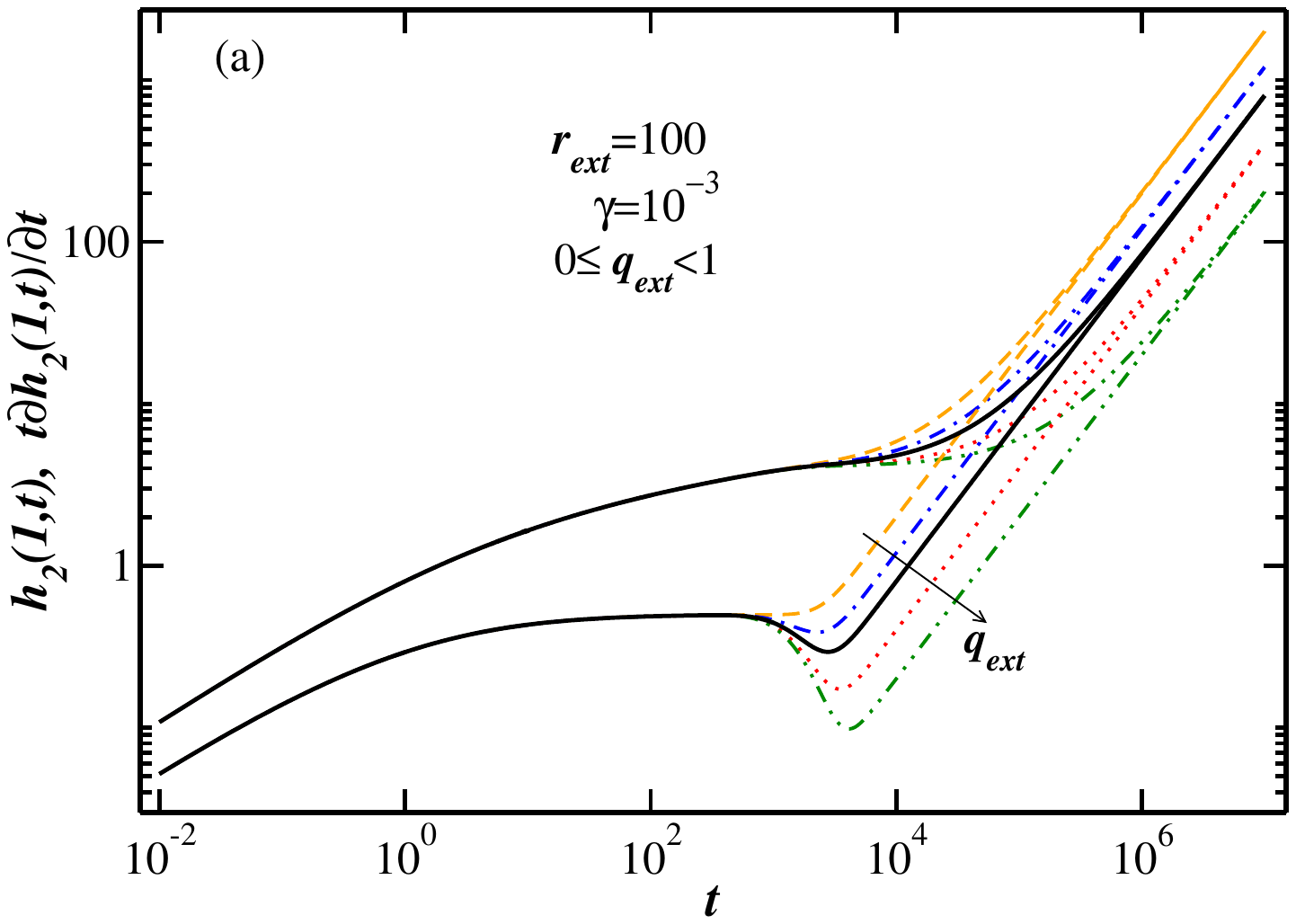}}
{\includegraphics[width=0.45\linewidth]{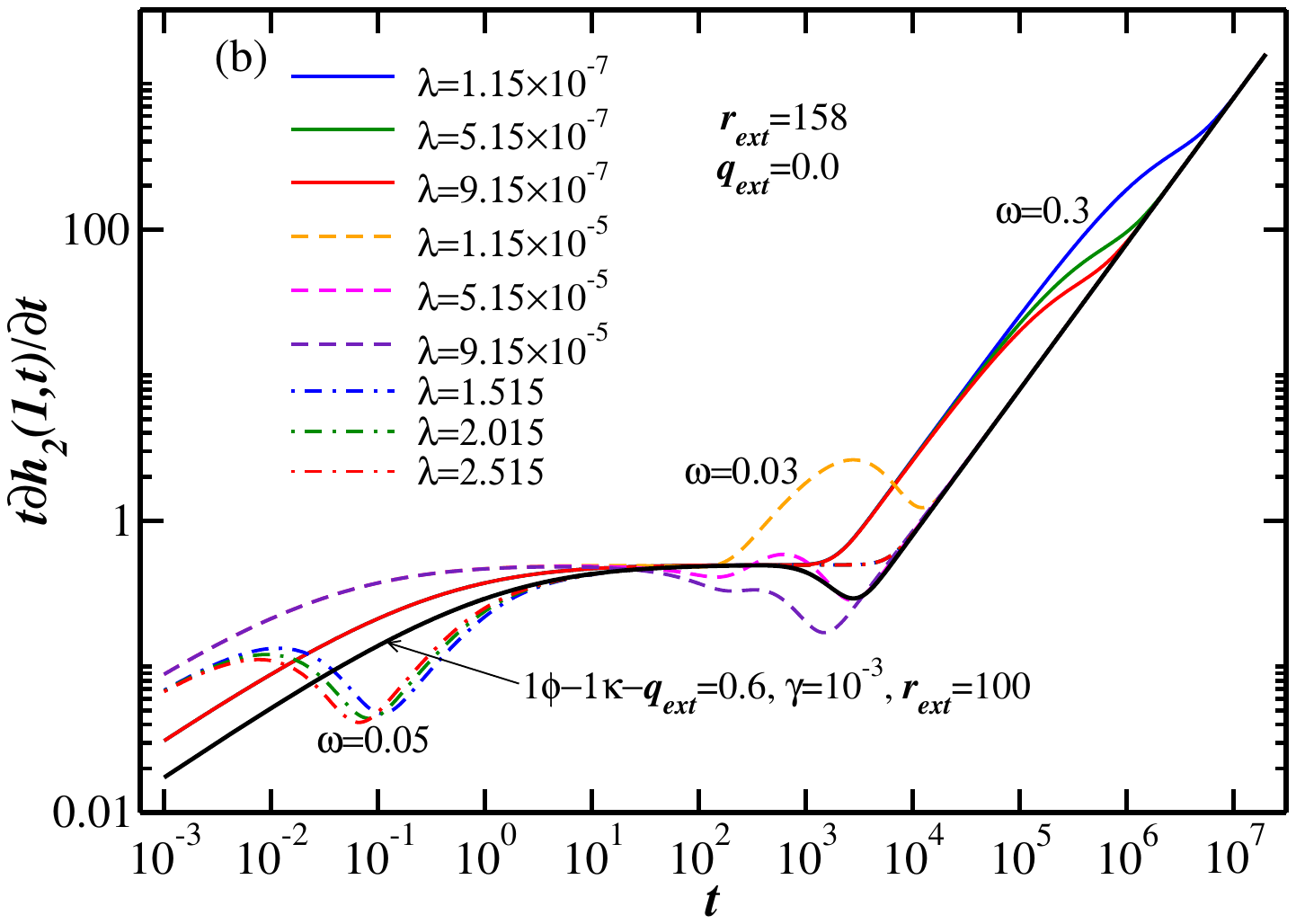}}\\
{\includegraphics[width=0.45\linewidth]{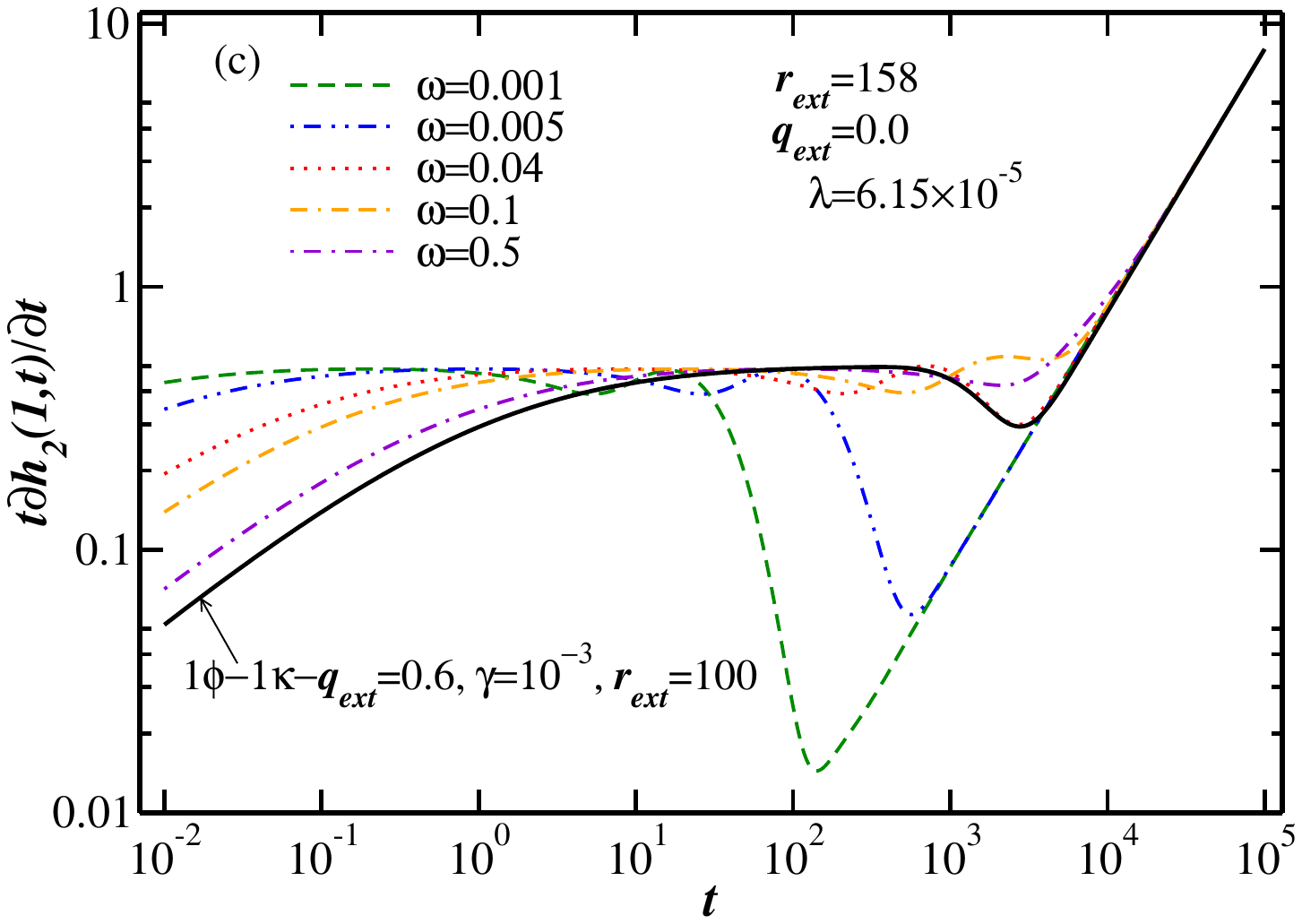}}
{\includegraphics[width=0.45\linewidth]{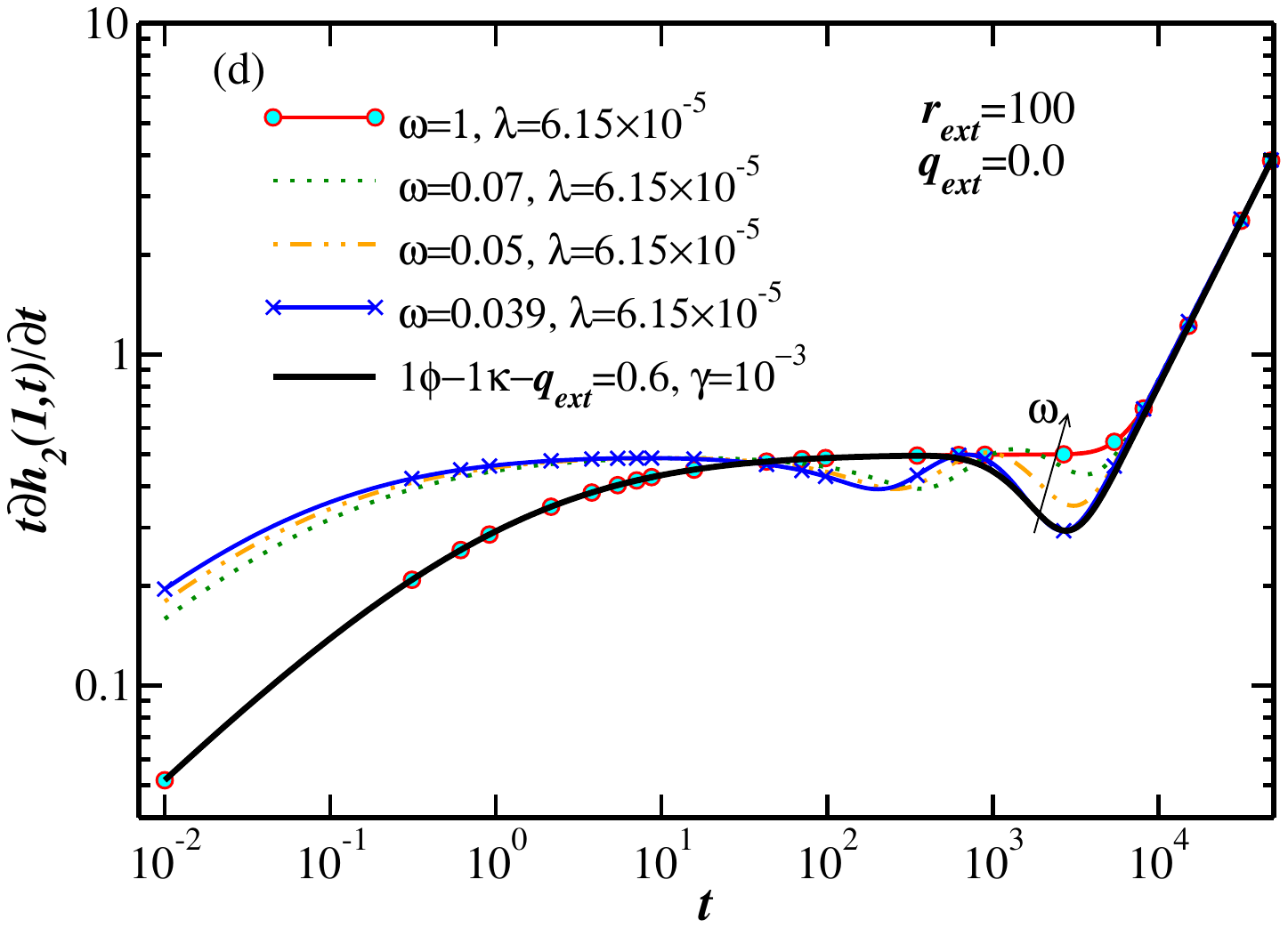}}
\caption{In (a) are presented graphs of the $1\phi$-$1\kappa$-$q_{\text{ext}}\neq 0$ model. The solid-black line remains as reference in (b), (c), and (d). The frames (b) and (c) show graphs of the drawdown derivative of the $2\phi$-$1\kappa$-$q_{\text{ext}}=0$ model in order to exhibit the curve behaviors and their similarities with the curve of the $1\phi$-$1\kappa$-$q_{\text{ext}}\neq 0$ model. In (d) it is shown that the best-fitting curve of the $2\phi$-$1\kappa$-$q_{\text{ext}}=0$ model (red line with circle symbol) is a curve of the $1\phi$-$1\kappa$-$q_{\text{ext}}=0$ model, since $\omega=1$ by fitting
}
\label{characteristic_behavior_equivalence_NN-BC}
\end{figure}

\subsection{Similarities between models with and without influx recharge}\label{Section_equivalence}

In this section, we indicate a criteria to know when a drawdown curve is related to a reservoir with or without influx recharge. The analysis presented is necessary because in Refs.~(\cite{doublet1995decline,del2014pressure}) is mentioned that there is a similarity between a model of fluid flow in a single-porosity medium with influx recharge ($1\phi$-$1\kappa$-$q_{\text{ext}}\neq 0$) and a model of fluid flow in a closed double-porosity medium ($2\phi$-$1\kappa$-$q_{\text{ext}}$=0). However, a detailed analysis using the Bourdet derivative is not carried out in these works in order to elucidate this statement, which is based on the fact that both models have characteristic curves with a minimum from a transition stage. Indeed, in this regard, note that in~\citet{wang2017transient} are shown curves of a $3\phi$-$1\kappa$-$q_{\text{ext}}\neq 0$ model that have a triple minimum in its drawdown derivative, i.e.~those curves may have equivalence with the curves of a fluid flow in a closed quad-porosity medium. The discussion may also occur because despite of the conceptual differences of the models being compared, their drawdown graphs are matched, as can be seen in Fig.~\ref{curves_with_or_without_recharge} (a) and (b). Therein, we compared the results from the following models: (a) $1\phi$-$1\kappa$-$q_{\text{ext}}\neq 0$ vs $2\phi$-$1\kappa$-$q_{\text{ext}}=0$ and (b) $2\phi$-$1\kappa$-$q_{\text{ext}}\neq 0$ vs $3\phi$-$1\kappa$-$q_{\text{ext}} = 0$. However, in the same figure is seen that the Bourdet derivative exhibits a clear difference in the transition period toward the outer boundary flow regime (stage with a unit slope at long times). Even so, the question remains whether it is possible to obtain equivalent curves from the results of $1\phi$-$1\kappa$-$q_{\text{ext}}\neq 0$ and $2\phi$-$1\kappa$-$q_{\text{ext}}=0$ models.

In order to elucidate the differences between $1\phi$-$1\kappa$-$q_{\text{ext}}\neq 0$ and $2\phi$-$1\kappa$-$q_{\text{ext}} = 0$ models, a comparison between both is shown in  Fig.~\ref{characteristic_behavior_equivalence_NN-BC}. This figure contains a characteristic drawdown curve (and/or its Bourdet derivative) of the $1\phi$-$1\kappa$-$q_{\text{ext}}\neq 0$ model (see solid-black lines) in order to remark the differences with curves of the $2\phi$-$1\kappa$-$q_{\text{ext}}=0$ model. Fig.~\ref{characteristic_behavior_equivalence_NN-BC} (a) shows drawdown curves with a monotonously increasing behavior, while the drawdown derivative has a minimum that increases its depth when the value of $q_{\text{ext}}$ increases ($0 \leq q_{\text{ext}} < 1$). In addition, a nonzero influx recharge leads to obtain a minimum located just before the stage dominated by the outer boundary effect. $2\phi$-$1\kappa$-$q_{\text{ext}}=0$ model also presents this latter stage, but its minimum, between the fractures and matrix flow regimes, is not necessarily located just before the beginning of the stage with a unit slope. This is exhibited in Fig.~\ref{characteristic_behavior_equivalence_NN-BC} (b) and (c), in which the systematic effect of the parameters $\lambda$ and $\omega$ on drawdown derivatives is observed for this model. In both frames there are derivative curves with two minimums, i.e.~these curves cannot be equivalent with the results in Fig.~\ref{characteristic_behavior_equivalence_NN-BC} (a), which clearly has a minimum. The rest of the curves in Fig.~\ref{characteristic_behavior_equivalence_NN-BC} (b) show behaviors which are dominated by the transition period between fractures and matrix flow regimes, which are present at short times when $\omega$=0.05 and at long times when $\omega$=0.3, see dashed-point and solid lines, respectively. Noticeably, these latter curves are very different of the ones in  Fig.~\ref{characteristic_behavior_equivalence_NN-BC} (a). Meanwhile, in Fig.~\ref{characteristic_behavior_equivalence_NN-BC} (c), $\omega$ is varied for the purpose of obtaining a curve with a single minimum and with the restriction that the transition period, between the fractures and matrix flow regimes, remains just at the beginning of the stage with a unit slope. As we can see, we can not eliminate one of the minimums by varying $\omega$. For this reason, we solved a least-square minimization, where the results of the $1\phi$-$1\kappa$-$q_{\text{ext}}\neq 0$ model are taken as the input data and the parametrized model is the $2\phi$-$1\kappa$-$q_{\text{ext}} = 0$. In this way, the curve-fitting problem leads to a value of $\omega$=1, i.e.~the best fit is a curve of the $1\phi$-$1\kappa$-$q_{\text{ext}}=0$ model [see red line with circles in Fig.~\ref{characteristic_behavior_equivalence_NN-BC} (d)]. In addition, we realize another fit by considering results  of the $1\phi$-$1\kappa$-$q_{\text{ext}}\neq 0$ model as input data, which only are taken from the minimum and the stage with unit slope. The best-fitting curve reproduces the input data, but outside of the fitting interval the curve has another minimum and at short times it is very different of the results of the $1\phi$-$1\kappa$-$q_{\text{ext}}=0$ model, see blue line with cross symbols in Fig.~\ref{characteristic_behavior_equivalence_NN-BC} (d). The rest of curves in Fig.~\ref{characteristic_behavior_equivalence_NN-BC} (d) are given in order to show the behavior by varying $\omega$; they have two minimums.

Therefore, according to our analysis it is impossible to obtain a quantitative equivalence between results of the $1\phi$-$1\kappa$-$q_{\text{ext}}\neq 0$ and $2\phi$-$1\kappa$-$q_{\text{ext}}=0$ models. We know that for $1\phi$-$1\kappa$-$q_{\text{ext}}\neq 0$ model, the only transition period is located at the beginning of the stage with unit slope. This stage has or has not a minimum when there is or there is not an influx recharge, respectively.  On the other hand, when there is a transition period that is not located at the beginning of the stage with unit slope, the porous medium is double-porosity, keeping out the possibility of quantitative equivalence with the $1\phi$-$1\kappa$-$q_{\text{ext}}\neq 0$ model. However, note that when the only transition period is located at the beginning of the stage with unit slope, but the transition has two minimums, the reservoir has associated a double-porosity medium without influx recharge. This latter statement can be used as a criterion to distinguish real-life data from a double-porosity reservoir with closed boundary. By contrast, we remark that a single-porosity reservoir with influx recharge has a minimum at the beginning of the stage with unit slope.

It is worth mentioning that by simplicity we compare $1\phi$-$1\kappa$-$q_{\text{ext}}\neq 0$ and $2\phi$-$1\kappa$-$q_{\text{ext}}=0$ models, on the understanding that similar conclusions must be obtained for other models, e.g.~$2\phi$-$1\kappa$-$q_{\text{ext}}\neq 0$ vs $3\phi$-$1\kappa$-$q_{\text{ext}}=0$.

\section{Conclusions}\label{Section_conclusions}

Using the joint Laplace-Hankel transform, we solved a fluid flow problem of interest in petroleum engineering and groundwater science. Our model describes the flow of a slightly compressible fluid in a double-porosity finite reservoir with Dirichlet-Dirichlet, Dirichlet-Neumann, Neumann-Dirichlet, and Neumann-Neumann boundary conditions. With this aim, the solution is divided in a stationary and a transient part. We validate the exact solution using the Stehfest method. In addition, when $\omega$=1 and $\lambda$=0, our formulas are reduced to those derived by other authors who solved the equivalent problem for a single-porosity model~(\cite{del2014pressure}).

We find that the Cinelli solution, Eq.~(\ref{inverse_h_2(k,s)d}), related to NN-BCs, is incomplete because it does not include time dependent terms, which impose the solution behavior at long time. For this reason, results of Eq.~(\ref{inverse_h_2(k,s)d}) always are stationary. The exact solution, in the limit of large $t$, is found by means of a series expansion in the Laplace space, thus, we identify the terms ignored by the Eq.~(\ref{inverse_h_2(k,s)d}), which are included in Eq.~(\ref{solWR_case4}). On the other hand, regarding the DD-BCs, DN-BCs, and ND-BCs cases, we observed a correct convergence of solutions (\ref{inverse_h_2(k,s)a}), (\ref{inverse_h_2(k,s)b}), and (\ref{inverse_h_2(k,s)c}), along the domain of solution, this except in the inner boundary where a zero value is obtained. Furthermore, these solutions are oscillatory and slowly converging. We find that the inhomogeneous BCs are reproduced using the closed formulas derived in this work, in addition, these simplified formulas speed up the convergence of the solutions compared to a direct use of the Cinelli relationships, Eqs.~(\ref{inverse_h_2(k,s)}).

Finally, the characteristic behaviors of solutions~(\ref{solWR_case1}), (\ref{solWR_case2}), (\ref{solWR_case3}), and (\ref{solWR_case4}), exhibit the different stages of flow at bottomhole: fractured-dominated, transitions-dominated, matrix-dominated, or recharge-dominated. The Bournet derivative, for NN-BCs case, shows a minimum during the transition periods between the fractures and matrix flow regimes and the dominated by the influx recharge flow regime. This information can be used to give a criterion about whether the reservoir has recharge at the outer boundary, at the same time that we can know the number of porosities associated with it.

\begin{acknowledgements}
L.X.~Vivas-Cruz thanks CONACYT (Mexico) for its financial support through a Ph.D. fellowship received during the realization of this work. A.~Gonz\'alez-Calder\'on and J.A.~Perera-Burgos acknowledge the support provided by CONACYT: C\'atedras CONACYT para j\'ovenes investigadores. The authors acknowledge to Felipe Pacheco-V\'azquez, Wilberth Herrera and Yarith del Angel for their comments and suggestions. The data can be found on Mendeley Datasets with the following doi: 10.17632/3837yfr46n.2.
\end{acknowledgements}

\section*{Conflict of interest}
The authors declare that they have no conflict of interest.

\renewcommand{\theequation}{A-\arabic{equation}}
\setcounter{equation}{0}

\section*{Appendix A: Exact solutions}

In this Appendix, we present additional details of the procedure developed in this work to solve the partial differential equation of the studied model.

The finite Hankel transform of a well-behaved function $h(r)$ is expressed as follows~(\cite{cinelli1965extension}):
\begin{equation}\label{CinelliFHT}
\tilde{h} (k_i) = H[h(r)] = \int_{a}^{b} rh(r)K(r,k_i,a)dr
\end{equation}
where $K(r,k_i,a)$ is a kernel that depends on the inner BC:
\begin{linenomath*}
\begin{subequations}
\begin{align}
\displaystyle K(r,k_i,a) &= \displaystyle J_0(k_i r)Y_0(k_ia) - J_0(k_i a)Y_0(k_i r), && \mbox{for DD and DN-BCs},\\
 \displaystyle &=J_1(k_i a)Y_0(k_i r) - J_0(k_i r)Y_1(k_ia), && \mbox{for ND and NN-BCs}.
 \end{align}
\end{subequations}
\end{linenomath*}
There is also an equivalent kernel that depends only on the outer boundary. Using the previous definition, the finite Hankel transform of the laplacian is given by
\begin{subequations}
\begin{align}
H\left[\frac{d^2h}{dr^2} + \frac{1}{r}\frac{dh}{dr}\right] =& \frac{2}{\pi}\frac{J_0(k_i a)}{J_0(k_i b)}h(b) - \frac{2}{\pi}h(a) - k_i^2\tilde{h}(k_i),\quad &&\mbox{for DD-BCs},\\
=& -\frac{2}{\pi k_i}\frac{J_0(k_i a)}{J_1(k_i b)}h'(b) - \frac{2}{\pi}h(a) - k_i^2\tilde{h}(k_i),\quad &&\mbox{for DN-BCs},\\
=& -\frac{2}{\pi k_i}h'(a) - \frac{2J_1(k_i a)}{\pi J_0(k_i b)}h(b) - k_i^2\tilde{h}(k_i),\quad &&\mbox{for ND-BCs},\\
=& \frac{2J_1(k_i a)}{\pi k_i J_1(k_i b)}h'(b) - \frac{2}{\pi k_i}h'(a) - k_i^2\tilde{h}(k_i),\quad &&\mbox{for NN-BCs}.
\end{align}
\end{subequations}
where $h'(a)$ is the derivative of $h(r)$ evaluated in $r=a$.
These latter expressions are used to obtain the inverse finite Hankel transform given in Eq.~\eqref{inverse_h_2(k,s)}.

\subsection*{Finite reservoir with Dirichlet-Dirichlet boundary conditions}
Using the JLHT to obtain the solution of a fluid flow in a finite reservoir with constant pressure in both boundaries is very simple. From Eq.~(\ref{h_2(k,s)a}) and Eqs.~(\ref{condition1}) and (\ref{condition3}), we get
\begin{linenomath*}
\begin{equation}\label{joint_HL_A}
\tilde{\widehat{h}}_{2}(k_i,s) = -\frac{1}{s} \mathcal{F}(k_i,s),
\end{equation}
\end{linenomath*}
where
\begin{linenomath*}
\begin{equation}
\mathcal{F}(k_i,s) = \displaystyle \frac{2}{ \pi [\eta(s) + k_i^2] },
\end{equation}
\end{linenomath*}
$\eta(s)$ being defined in Eq.~(\ref{eta(s)})

Taking the inverse Laplace transform of Eq.~(\ref{joint_HL_A}) leads to the following expression in Hankel space:
\begin{linenomath*}
\begin{equation}\label{h2(k,t)_2Pmodel}
\tilde{h}_{2}(k_i,t) = -\frac{2}{\pi k_i^2} + \tilde{g}(k_i,t),
\end{equation}
\end{linenomath*}
where
\begin{linenomath*}
\begin{eqnarray}
\tilde{g}(k_i,t) &=& \frac{1}{\nu}\big[(\nu + \varrho)\exp\left\lbrace -\displaystyle \frac{\xi + \nu}{2\psi} t\right\rbrace + (\nu - \varrho) \displaystyle \exp\left\lbrace-\frac{\xi-\nu}{2\psi}t\right\rbrace \big],\label{g_D1}
\end{eqnarray}
\end{linenomath*}
and
\begin{linenomath*}
\begin{equation}\label{functions_nu_xi_varrho}
\begin{split}
\psi & = \omega(\omega-1),\\
\xi & = -\lambda + k_i^2(\omega - 1),\\
\varrho & = \lambda + k_i^2(\omega - 1),\\
\nu & = \sqrt{\xi^2 + 4k_i^2\lambda\psi}.\\
\end{split}
\end{equation}
\end{linenomath*}
Because the term $-2/\pi k_i^2$ is not related to time, its finite inverse Hankel transform is the stationary solution of model~(\ref{2Pmodel}). Therefore, we replace the inversion of $-2/\pi k_i^2$ by the stationary solution of the Laplace equation, $\nabla^2h_2=0$, in cylindrical coordinates. Namely, the following equality holds true:
\begin{linenomath*}
\begin{eqnarray}
1 - \frac{\log(r)}{\log(r_{\text{ext}})} &=& -\pi \sum_{i=1}^{\infty} \frac{\mathcal{I}_{0,0}(k_i,r,1) J_0^2(r_{\text{ext}}k_i)}{J_0^2(k_i) - J_0^2(r_{\text{ext}}k_i)}.\label{term1}
\end{eqnarray}
\end{linenomath*}
Substituting Eq.~(\ref{h2(k,t)_2Pmodel}) into Eq.~(\ref{inverse_h_2(k,s)}) for the DD-BCs case, and simplifying with the previous closed formula, the exact solution is
\begin{linenomath*}
\begin{equation}\label{solWR_case1Ap}
h_{2}(r,t) = 1 - \frac{\log(r)}{\log(r_{\text{ext}})} + \frac{\pi}{2} \sum_{i=1}^{\infty} \frac{ \tilde{g}(k_i,t) \mathcal{I}_{0,0}(k_i,r,1) J_0^2(r_{\text{ext}}k_i)}{J_0^2(k_i) - J_0^2(r_{\text{ext}}k_i)},
\end{equation}
\end{linenomath*}
where the $k_i$'$s$ are the positive roots of $\mathcal{I}_{0,0}(k_i,1,r_{\text{ext}})=0$.

The flux is obtained by substituting Eq.~(\ref{solWR_case1}) into Eq.~(\ref{flux_Darcylaw}). This is accomplished as follows:
\begin{linenomath*}
\begin{equation}\label{solFlux_case1Ap}
j_{2}(t) = \frac{1}{\log(r_{\text{ext}})} - \frac{\pi}{2} \sum_{i=1}^{\infty} \frac{k_i \tilde{g}(k_i,t) \mathcal{I}_{0,1}(k_i,1,1) J_0^2(r_{\text{ext}}k_i)}{J_0^2(k_i) - J_0^2(r_{\text{ext}}k_i)}.
\end{equation}
\end{linenomath*}
Eqs.~(\ref{solWR_case1}) and (\ref{solFlux_case1})  recover the formulas in~\citet{muskat1934flow} when the limit of single-porosity medium is taken.

\subsection*{Finite reservoir with Dirichlet-Neumann boundary conditions}
Fluid flow in a reservoir with constant pressure at the bottomhole and influx recharge $f(t)$ at the outer boundary is considered. The influx function $f(t)$ is given in Eq.~(\ref{constant_rate}). Substituting Eqs.~(\ref{condition1}) and (\ref{condition4}) into Eq.~(\ref{h_2(k,s)b}) leads to the following formula:

\begin{linenomath*}
\begin{equation}\label{joint_HL_B}
\tilde{\widehat{h}}_{2}(k_i,s) = \left( \frac{q_{\text{ext}}J_0(k_i)}{r_{\text{ext}} k_i J_1(r_{\text{ext}}k_i)} \frac{1}{(\gamma s^2 + s)} - \frac{1}{s}\right) \mathcal{F}(k_i,s).
\end{equation}
\end{linenomath*}

Taking the inverse Laplace transform of Eq.~(\ref{joint_HL_B}), we obtain
\begin{linenomath*}
\begin{equation}\label{h2(k,t)_caseB}
\begin{split}
\tilde{h}_{2}(k_i,t) =& \displaystyle \frac{2}{\pi k_i^2}\left[ \frac{q_{\text{ext}}}{r_{\text{ext}}}\frac{J_0(k_i)}{k_i J_1(r_{\text{ext}}k_i)} - 1 \right] + \tilde{g}(k_i,t),
\end{split}
\end{equation}
\end{linenomath*}
where
\begin{linenomath*}
\begin{eqnarray}
\tilde{g}(k_i,t) &=& \frac{q_{\text{ext}} J_0(k_i)}{\pi k_i^3 r_{\text{ext}} J_1(r_{\text{ext}}k_i) \vartheta \nu}\exp\left\lbrace -\frac{\xi+\nu}{2\psi}t \right\rbrace \Big[ \Big(-2k_i^2 \gamma(\gamma\lambda + \omega-1) \exp \left\lbrace (\frac{\xi+\nu}{2\psi} - \frac{1}{\gamma}) t \right\rbrace + \nonumber\\[5pt]
&&(\psi+\lambda\gamma)(\exp \left\lbrace \frac{\nu}{\psi}t \right\rbrace +1)\Big) \nu +(\lambda\gamma (\xi-2k_i^2\psi) - \psi \varrho)\left(\exp \left\lbrace \frac{\nu}{\psi}t \right\rbrace -1 \right)\Big]-\label{wD_caseB}\\[5pt]
&&\frac{1}{\pi k_i^2 \nu}\exp\left\lbrace - \frac{\xi + \nu}{2\psi} t\right\rbrace \left[ \left(\exp\left\lbrace\frac{\nu}{\psi}t\right\rbrace - 1 \right)\varrho - \left( \exp\left\lbrace \frac{\nu}{\psi}t\right\rbrace + 1 \right)\nu \right],\nonumber
\end{eqnarray}
\end{linenomath*}
and $\xi$, $\psi$, $\varrho$ and $\nu$, are given in Eq.~(\ref{functions_nu_xi_varrho}). In addition, we have
\begin{linenomath*}
\begin{eqnarray}
\vartheta = \gamma \xi - \psi + k_i^2 \gamma^2 \lambda.
\end{eqnarray}
\end{linenomath*}
Similar to DD-BCs case, Eq.~(\ref{h2(k,t)_caseB}) has terms not related to time, whose inverse finite Hankel transform is the stationary solution of model~(\ref{2Pmodel}). Therefore, this inverse is equal to the solution of Laplace equation with DN-BCs:
\begin{linenomath*}
\begin{eqnarray}
1 - q_{\text{ext}}\log (r) &=& \pi\sum_{i=1}^{\infty}\Big( \frac{q_{\text{ext}}}{r_{\text{ext}}}\frac{J_0(k_i)}{k_i J_1(r_{\text{ext}}k_i)} - 1 \Big) \frac{\mathcal{I}_{0,0}(k_i,r,1) J_1^2(r_{\text{ext}}k_i)}{J_0^2(k_i) - J_1^2(r_{\text{ext}}k_i)}.\label{term2}
\end{eqnarray}
\end{linenomath*}
Substituting Eq.~(\ref{h2(k,t)_caseB}) into Eq.~(\ref{inverse_h_2(k,s)b}), and simplifying with the closed formula~(\ref{term2}), the exact solution is
\begin{linenomath*}
\begin{equation}\label{solWR_case2Ap}
h_{2}(r,t) = 1 - q_{\text{ext}}\log(r) + \frac{\pi^2}{2} \sum_{i=1}^{\infty} \frac{ k_i^2 \tilde{g}(k_i,t) \mathcal{I}_{0,0}(k_i,r,1) J_1^2(r_{\text{ext}}k_i)}{J_0^2(k_i) - J_1^2(r_{\text{ext}}k_i)},
\end{equation}
\end{linenomath*}
where the $k_i$'$s$ are the positive roots of $\mathcal{I}_{1,0}(k_i,r_{\text{ext}},1)=0$.

In addition, from Eqs.~(\ref{flux_Darcylaw}) and (\ref{solWR_case2}), the flux at the bottomhole can be written as
\begin{linenomath*}
\begin{equation}\label{Radialflux_case2Ap}
j_{2}(t) = q_{\text{ext}} - \frac{\pi^2}{2} \sum_{i=1}^{\infty} \frac{ k_i^3 \tilde{g}(k_i,t) \mathcal{I}_{0,1}(k_i,1,1) J_1^2(r_{\text{ext}}k_i)}{J_0^2(k_i) - J_1^2(r_{\text{ext}}k_i)}.
\end{equation}
\end{linenomath*}
Eqs.~(\ref{solWR_case2}) and (\ref{Radialflux_case2}) recover the formulas in~\citet{muskat1934flow} and \citet{hurst1934unsteady} when the limit of single-porosity is taken.

\subsection*{Finite reservoir with Neumann-Dirichlet boundary conditions}
Next, assuming a constant terminal rate at the bottomhole and constant pressure at the outer boundary of a finite reservoir, the pressure is found; i.e.~the ND-BCs case is solved. Substituting Eqs.~(\ref{condition2}) and (\ref{condition3}) in Eq.~(\ref{h_2(k,s)c}), we obtain that
\begin{linenomath*}
\begin{equation}\label{joint_HL_C}
\tilde{\widehat{h}}_{2}(k_i,s) = \frac{1}{k_i s}\mathcal{F}(k_i,s).
\end{equation}
\end{linenomath*}
Taking the inverse Laplace transform of Eq.~(\ref{joint_HL_C}) leads to
\begin{linenomath*}
\begin{equation}\label{h2(k,t)_2Pmodel_case3}
\tilde{h}_{2}(k_i,t) = \frac{2}{\pi k_i^3} + \tilde{g}(k_i,t),
\end{equation}
\end{linenomath*}
where
\begin{linenomath*}
\begin{eqnarray}
\tilde{g}(k_i,t) &=& \frac{1}{\nu}\Big[(\varrho-\nu)\exp\left\lbrace \displaystyle \frac{\nu-\xi}{2\psi} t\right\rbrace - (\varrho+\nu) \displaystyle \exp\left\lbrace -\frac{\xi+\nu}{2\psi}t\right\rbrace \Big].\label{g_D}
\end{eqnarray}
\end{linenomath*}
From the inverse of terms not related to time in Eq.~(\ref{h2(k,t)_2Pmodel_case3}) and from the time-independent solution of the model~(\ref{2Pmodel}), we have
\begin{linenomath*}
\begin{eqnarray}
\log\left( \frac{r_{\text{ext}}}{r} \right) &=& \pi \sum_{i=1}^{\infty} \frac{ \mathcal{I}_{1,0}(k_i,1,r) J_0^2(r_{\text{ext}}k_i)}{k_i[J_1^2(k_i) - J_0^2(r_{\text{ext}}k_i)]}.\label{term3}
\end{eqnarray}
\end{linenomath*}
Replacing Eq.~(\ref{h2(k,t)_2Pmodel_case3}) in Eq.~(\ref{inverse_h_2(k,s)c}), and using Eq.~({\ref{term3}), the exact solution is given by
\begin{linenomath*}
\begin{equation}\label{solWR_case3Ap}
h_{2}(r,t) = \log\left( \frac{r_{\text{ext}}}{r} \right) + \frac{\pi}{2} \sum_{i=1}^{\infty} \frac{ \tilde{g}(k_i,t) \mathcal{I}_{1,0}(k_i,1,r) J_0^2(r_{\text{ext}}k_i)}{k_i[J_1^2(k_i) - J_0^2(r_{\text{ext}}k_i)]},
\end{equation}
\end{linenomath*}
where the $k_i$'$s$ are the positive roots of $\mathcal{I}_{1,0}(k_i,1,r_{\text{ext}})=0$.

Eq.~(\ref{solWR_case3}) recovers the formulas in~\citet{matthews1967pressure} when the limit of single-porosity is taken.

\subsection*{Finite reservoir with Neumann-Neumann boundary conditions}
The pressure of a fluid in a double-porosity reservoir with constant terminal rate and influx recharge $f(t)$ at the outer boundary is given. Replacing Eqs.~(\ref{condition2}) and (\ref{condition4}) in Eq.~(\ref{h_2(k,s)d}), we obtain
\begin{linenomath*}
\begin{equation}\label{joint_HL_case4_2}
\begin{array}{lll}
\displaystyle\tilde{\widehat{h}}_{2}(k_i,s) &=& \displaystyle\frac{\mathcal{F}(k_i,s)}{k_i s} - q_{\text{ext}} \frac{ J_1(k_i)\mathcal{F}(k_i,s)}{r_{\text{ext}} k_i (\gamma s^2 + s) J_1(r_{\text{ext}}k_i)}, \\[12pt]
&=& \displaystyle \tilde{\widehat{h}}_{2,\mbox{nf}}+\tilde{\widehat{h}}_{2,\mbox{f}}.
\end{array}
\end{equation}
\end{linenomath*}
The no-flow term $\tilde{\widehat{h}}_{2,\mbox{nf}}$ has the following inverse Laplace transform:
\begin{linenomath*}
\begin{equation}\label{h2(k,t)_2Pmodel_case4}
\displaystyle \tilde{h}_{2,\mbox{nf}}(k_i,t) = \frac{2}{\pi k_i^3} + \tilde{g}(k_i,t),
\end{equation}
\end{linenomath*}
where $\tilde{g}(k_i,t)$ is equal to the expression in Eq.~(\ref{g_D}), with the difference that the values of $k_i$'$s$  comes from the positive roots of $\mathcal{I}_{1,1}(k_i,1,r_{\text{ext}})$.

The second term in the RHS of Eq.~(\ref{joint_HL_case4_2}) is an influx term, whose inverse Laplace transform is found using the convolution theorem. This inverse can be written as
\begin{linenomath*}
\begin{eqnarray}
\tilde{h}_{2,\mbox{f}} & = & \displaystyle \frac{1}{r_{\text{ext}}}\tilde{\mathcal{G}}(k_i)\left( \mathcal{L}^{-1}\left\lbrace \widehat{f}(s)\right\rbrace \ast \mathcal{L}^{-1}\left\lbrace \tilde{\widehat{u}}(k_i,s) \right\rbrace \right), \nonumber \\
&=& \displaystyle \frac{\tilde{\mathcal{G}}(k_i)}{r_{\text{ext}}} \int_0^t f(t-\zeta)\tilde{u}(k_i,\zeta)d\zeta,\nonumber\\
&=& \displaystyle -\frac{q_{\text{ext}}}{r_{\text{ext}}}\tilde{\mathcal{G}}(k_i)\left[ \tilde{Q}_1(k_i,t) - \tilde{Q}_2(k_i,t) \right],\label{Q1_Q2}
\end{eqnarray}
\end{linenomath*}
in which we use the following inverse Laplace transform
\begin{linenomath*}
\begin{equation}\label{u(k_i,t)}
\tilde{u}(k_i,t) = \frac{1}{2\omega \nu}\exp\left\lbrace -\frac{\xi+\nu}{2\psi}t\right\rbrace\Big[ \left( 1 - \exp\left\lbrace \frac{\nu}{\psi}t\right\rbrace \right) (\xi + 2\omega\lambda) + \left( 1 + \exp\left\lbrace \frac{\nu}{\psi}t\right\rbrace \right)\nu  \Big].
\end{equation}
\end{linenomath*}
In previous equations, we have
\begin{linenomath*}
\begin{equation}
\begin{array}{lll}
\displaystyle \tilde{\mathcal{G}}(k_i)&=&\displaystyle \frac{2J_1(k_i)}{\pi k_i J_1(r_{\text{ext}}k_i)}\\[10pt]
\displaystyle \tilde{\widehat{u}}(k_i,s) &=& \displaystyle \frac{1}{\eta(s) + k_i^2},\\[12pt]
\displaystyle \tilde{Q}_1(k_i,t) &=& \displaystyle \int_0^t \tilde{u}(k_i,\zeta)d\zeta\\
& = & \displaystyle \frac{1}{2 k_i^2 \nu}\exp\left\lbrace -\displaystyle \frac{(\xi + \nu)}{2\psi} t\right\rbrace \left[ \displaystyle \left(\exp\left\lbrace\frac{\nu}{\psi}t\right\rbrace - 1 \right)\varrho + \left( 2\exp\left\lbrace \frac{(\xi + \nu)}{2\psi}t\right\rbrace - \exp\left\lbrace \frac{\nu}{\psi}t\right\rbrace -1 \right)\nu \right],\\[12pt]
\displaystyle \tilde{Q}_2(k_i,t) &=& \displaystyle \int_0^t \mbox{e}^{-(t-\zeta)/\gamma}\tilde{u}(k_i,\zeta)d\zeta\\[12pt]
& = & \displaystyle \frac{\gamma (\omega-1)}{\nu}\mbox{e}^{-t/\gamma} \Big[ \frac{\mathcal{A}+\nu}{\mathcal{B}-2\psi} - \frac{\mathcal{A}-\nu}{\mathcal{C}-2\psi} - \mbox{e}^{-(\mathcal{B}-2\psi)/(2\gamma\psi)t} \left( \frac{\mbox{e}^{(\nu/\psi) t}(\mathcal{A}-\nu)}{\mathcal{C}+2\psi} + \frac{\mathcal{A}+\nu}{\mathcal{B}-2\psi} \right) \Big],
\end{array}
\end{equation}
\end{linenomath*}
where $\xi$, $\nu$ and $\psi$, are given in Eqs.~(\ref{functions_nu_xi_varrho}), $\mathcal{A} = \xi+2\omega \lambda$, $\mathcal{B} = (\nu+\xi)\gamma$ and $\mathcal{C} = (\nu-\xi)\gamma$.

Therefore, the inverse Laplace transform of Eq.~(\ref{joint_HL_case4_2}) is the sum of Eqs.~(\ref{h2(k,t)_2Pmodel_case4}) and (\ref{Q1_Q2}):
\begin{linenomath*}
\begin{eqnarray}\label{Eqs.40-45}
\tilde{h}_{2}(k_i,t) = \tilde{h}_{2,\mbox{nf}} +\tilde{h}_{2,\mbox{f}}.
\end{eqnarray}
\end{linenomath*}

Since in the NN-BCs case the Laplace equation has no solution, the asymptotic solution for long time is found by means of a series expansion of the solution in Laplace space, i.e.~the expansion of Eq.~(\ref{WR_h2D_caseD}) about $s=0$~(\cite{van1949application,prats1986interpretation}) is developed. Thus, the behavior for long time of the no-flow term was given by~\citet{van1949application}:
\begin{linenomath*}
\begin{equation}\label{Eq.VII-13_vanEverd}
h_{2,\text{n-f}}(r,t) = \frac{2}{r_{\text{ext}}^2-1}\left( \frac{r^2}{4}+t \right) - \frac{r_{\text{ext}}^2}{r_{\text{ext}}^2-1}\log (r) - \frac{3r_{\text{ext}}^4 - 4r_{\text{ext}}^4 \log (r_{\text{ext}})-2r_{\text{ext}}^2-1}{4(r_{\text{ext}}^2-1)^2},
\end{equation}
\end{linenomath*}
while for the influx term was found by~\citet{del2014pressure}:
\begin{linenomath*}
\begin{equation}\label{q-t_Large}
h_{2,\text{f}}(r,t) = q_{\text{ext}} \Big[ \frac{2}{r_{\text{ext}}^2-1}\left( \frac{r^2}{4} + t \right) - \frac{\log (r)}{r_{\text{ext}}^2-1}- \frac{r_{\text{ext}}^4 + 2r_{\text{ext}}^2 - 4r_{\text{ext}}^2 \log (r_{\text{ext}}) -3}{4(r_{\text{ext}}^2-1)^2} \Big] + \frac{2 \gamma q_{\text{ext}}(1 - \mbox{e}^{-t/\gamma})}{r_{\text{ext}}^2-1}.
\end{equation}
\end{linenomath*}
Equalities (\ref{Eq.VII-13_vanEverd}) and (\ref{q-t_Large}) come from studies of fluid flow in a single-porosity medium. However, they can be used in the solution of the double-porosity model, since at long time, the fluid behavior resembles that of a fluid in an homogeneous reservoir. Mathematically,
\begin{linenomath*}
\begin{equation}
\eta(s)\rightarrow s \mbox{ as } s\rightarrow 0,
\end{equation}
\end{linenomath*}
which is the $\eta(s)$ of a single-porosity medium. Accordingly, we can use the Eqs.~(\ref{Eq.VII-13_vanEverd}) and (\ref{q-t_Large}) in the model of double-porosity.

Eqs.~(\ref{h2(k,t)_2Pmodel_case3}) and (\ref{Q1_Q2}) include the terms $2/\pi k_i^3$ and $-q_{\text{ext}}\tilde{\mathcal{G}}(k_i)/(r_{\text{ext}}k_i^2)$, respectively, whose inverse finite Hankel transform is time-independent. This implies that the time-independent terms of Eqs.~(\ref{Eq.VII-13_vanEverd}) and (\ref{q-t_Large}) must equal the inverse finite Hankel transform of $2/\pi k_i^3$ and $-q_{\text{ext}}\tilde{\mathcal{G}}(k_i)/(r_{\text{ext}}k_i^2)$, respectively. Therefore,
\begin{linenomath*}
\begin{equation}\label{Term_NN_noflow}
\pi \sum_{i=1}^{\infty} \frac{\mathcal{I}_{1,0}(k_i,1,r) J_1^2(r_{\text{ext}}k_i)}{k_i [J_1^2(k_i) - J_1^2(r_{\text{ext}}k_i)]} = \frac{r^2}{2(r_{\text{ext}}^2-1)} - \frac{r_{\text{ext}}^2}{r_{\text{ext}}^2-1}\log (r) - \frac{3r_{\text{ext}}^4 - 4r_{\text{ext}}^4 \log (r_{\text{ext}})-2r_{\text{ext}}^2-1}{4(r_{\text{ext}}^2-1)^2},
\end{equation}
\end{linenomath*}
and
\begin{linenomath*}
\begin{eqnarray}
-q_{\text{ext}} \pi \sum_{i=1}^{\infty} \frac{J_1(k_i)\mathcal{I}_{1,0}(k_i,1,r) J_1(r_{\text{ext}}k_i)}{r_{\text{ext}}k_i [J_1^2(k_i) - J_1^2(r_{\text{ext}}k_i)]} &=& -q_{\text{ext}} \left[ \frac{r^2}{2(r_{\text{ext}}^2-1)} - \frac{\log (r)}{r_{\text{ext}}^2-1} - \frac{r_{\text{ext}}^4 + 2r_{\text{ext}}^2 - 4r_{\text{ext}}^2 \log (r_{\text{ext}}) -3}{4(r_{\text{ext}}^2-1)^2} \right] \nonumber\\
&&+ \frac{2 \gamma q_{\text{ext}}}{r_{\text{ext}}^2-1} \label{Term_NN_influx}.
\end{eqnarray}
\end{linenomath*}
Substituting Eq.~(\ref{Eqs.40-45}) into Eq.~(\ref{inverse_h_2(k,s)d}) and using the previous closed formulas, the pressure is
\begin{linenomath*}
\begin{equation}\label{solWR_case4Ap}
\begin{split}
h_{2}(r,t) =& \frac{\pi^2}{2} \sum_{i=1}^{\infty} \frac{k_i^2 \tilde{g}(k_i,t) \mathcal{I}_{1,0}(k_i,1,r) J_1^2(r_{\text{ext}}k_i)}{J_1^2(k_i) - J_1^2(r_{\text{ext}}k_i)}+ \frac{2}{r_{\text{ext}}^2-1}\left( \frac{r^2}{4}+t \right) - \frac{r_{\text{ext}}^2}{r_{\text{ext}}^2-1}\log (r)\\[3pt]
& - \frac{3r_{\text{ext}}^4 - 4r_{\text{ext}}^4 \log(r_{\text{ext}}) - 2r_{\text{ext}}^2 - 1}{4(r_{\text{ext}}^2-1)^2} + q_{\text{ext}} \Big[ \frac{2}{r_{\text{ext}}^2-1}\left( \frac{r^2}{4} + t \right) - \frac{\log (r)}{r_{\text{ext}}^2-1}\\[3pt]
& - \frac{r_{\text{ext}}^4 + 2r_{\text{ext}}^2 - 4r_{\text{ext}}^2 \log (r_{\text{ext}}) -3}{4(r_{\text{ext}}^2-1)^2} \Big] + \frac{2 \gamma q_{\text{ext}}(1 - \mbox{e}^{-t/\gamma})}{r_{\text{ext}}^2-1},
\end{split}
\end{equation}
\end{linenomath*}
where
\begin{linenomath*}
\begin{equation}\label{x*}
\begin{array}{lll}
\displaystyle \tilde{g}(k_i,t) &=& \displaystyle \tilde{\chi}(k_i,t) + (q_{\text{ext}}/r_{\text{ext}})\tilde{\mathcal{G}}(k_i)\big[\tilde{Q}_2(k_i,t)-\tilde{R}_1(k_i,t)\big], \\[12pt]
\displaystyle \tilde{\chi}(k_i,t) &=& \displaystyle \frac{1}{\pi k_i^3 \nu}\exp\left\lbrace -\displaystyle \frac{(\xi + \nu)}{2\psi} t\right\rbrace \left[ \displaystyle \left(\exp\left\lbrace\frac{\nu}{\psi}t\right\rbrace - 1 \right)\varrho - \left( \exp\left\lbrace \frac{\nu}{\psi}t\right\rbrace + 1 \right)\nu \right],\\[12pt]
\displaystyle \tilde{R}_1(k_i,t) &=& \displaystyle \frac{\pi }{2}k_i\tilde{\chi}(k_i,t),
\end{array}
\end{equation}
\end{linenomath*}
and $k_i$'$s$ are the positive roots of $\mathcal{I}_{1,1}(k_i,1,r_{\text{ext}})$. Note that to find $\tilde{g}(k_i,t)$, we write $\tilde{Q}_1(k_i,t) = 1/k_i^2 + \tilde{R}_1(k_i,t)$.

The time dependent terms in Eqs.~(\ref{Eq.VII-13_vanEverd}) and (\ref{q-t_Large}) are included in Eq.~(\ref{solWR_case4Ap}) in order to describe the long time fluid behavior. It is worth mentioning that these terms are omitted in the Cinelli solution~\cite{cinelli1965extension} of the NN-BCs case. Eq.~(\ref{solWR_case4Ap}) recovers the results in~\citet{muskat1934flow}, \citet{matthews1967pressure}, and \citet{del2014pressure}, when the limit of single-porosity is taken.

\renewcommand{\theequation}{B-\arabic{equation}}
\setcounter{equation}{0}
\section*{Appendix B: solutions in Laplace space for the different cases of study}
\renewcommand{\arraystretch}{0}

This Appendix presents the exact analytical results for pressure and flux, in the Laplace space, for our study model. Namely, pressure formulas in Table~\ref{TableAppendix} are the solutions of Eq.~(\ref{Laplace_2Dmodel_Bessel}) with different combinations of the BCs in (\ref{condition1})-(\ref{condition5}), while flux formulas come from replacing these relationships in Eq.~(\ref{flux_Darcylaw}).

\makeatletter
\setcounter{table}{0}
\renewcommand{\thetable}{A\@arabic\c@table}
\makeatother

\begin{table}[t]
\small\addtolength{\tabcolsep}{-10pt}
\caption{\label{TableAppendix}Exact solutions in Laplace space of the fluid flow equation~(\ref{2Pmodel}) and their different BCs cases.}
\begin{tabular}{p{1.3cm}p{8.5cm}p{8cm}}
\toprule[0.7pt]
\begin{center}Case\end{center} & \begin{center}Pressure$^*$\end{center} & \begin{center}Flux$^*$\end{center}\\[0.5pt]
\cmidrule[0.7pt](r){1-3} \\[5pt]
\begin{center}
DD-BC
\end{center}

&
\begin{eqnarray}
\widehat{h}_{2}(r,s) =\frac{\Psi_{0,0}\left(r_{\text{ext}},r,\sqrt{\eta(s)}\right)}{s \Psi_{0,0}\left(r_{\text{ext}},1,\sqrt{\eta(s)}\right)}\label{WR_h2D_caseA}
\end{eqnarray}

&
\begin{eqnarray}
\widehat{j}_{2}(s) =\frac{\sqrt{\eta(s)}\Psi_{0,1}\left(r_{\text{ext}},1,\sqrt{\eta(s)}\right)}{s \Psi_{0,0}\left(1,r_{\text{ext}},\sqrt{\eta(s)}\right)}\label{WR_j2D_caseA}
\end{eqnarray}

\\[5pt]
\begin{center}
DN-BC
\end{center}

&

\begin{eqnarray}\label{WR_h2D_caseB-1}
\begin{split}
\widehat{h}_{2}(r,s) &= \frac{\Psi_{1,0}\left(r_{\text{ext}},r,\sqrt{\eta(s)}\right)}{s \Psi_{0,1}\left(1,r_{\text{ext}},\sqrt{\eta(s)}\right)} + \\
&\hspace{0.3cm}\frac{\widehat{f}(s)}{r_{\text{ext}}}\frac{\Psi_{0,0}\left(1,r,\sqrt{\eta(s)}\right)}{\sqrt{\eta(s)} \Psi_{0,1}\left(1,r_{\text{ext}},\sqrt{\eta(s)}\right)}
\end{split}
\end{eqnarray}

&
\begin{eqnarray}\label{WR_h2D_caseB-2}
\begin{split}
\widehat{j}_{2}(s) &= \frac{\sqrt{\eta(s)}\Psi_{1,1}\left(1,r_{\text{ext}},\sqrt{\eta(s)}\right)}{s \Psi_{0,1}\left(1,r_{\text{ext}},\sqrt{\eta(s)}\right)} -\\
&\hspace{0.4cm}\frac{\widehat{f}(s)}{r_{\text{ext}}}\frac{\Psi_{0,1}\left(1,1,\sqrt{\eta(s)}\right)}{ \Psi_{0,1}\left(1,r_{\text{ext}},\sqrt{\eta(s)}\right)}
\end{split}
\end{eqnarray}
\\[5pt]

\begin{center}
ND-BC
\end{center}

&
\begin{equation}
\begin{split}
\widehat{h}_{2}(r,s) =\frac{\Psi_{0,0}\left(r,r_{\text{ext}},\sqrt{\eta(s)}\right)}{s \sqrt{\eta(s)}\Psi_{0,1}\left(r_{\text{ext}},1,\sqrt{\eta(s)}\right)}
\end{split}
\label{WR_h2D_caseC}
\end{equation}

&

\begin{eqnarray}\label{flux_Laplace_caseC}
\begin{split}
\widehat{j}_{2}(s) = \frac{1}{s}
\end{split}
\end{eqnarray}

\\[5pt]
\begin{center}
NN-BC
\end{center}

&
\begin{eqnarray}\label{WR_h2D_caseD}
\begin{split}
\widehat{h}_{2}(r,s) &= \frac{\Psi_{0,1}\left(r,r_{\text{ext}},\sqrt{\eta(s)}\right)}{s \sqrt{\eta(s)}\Psi_{1,1}\left(1,r_{\text{ext}},\sqrt{\eta(s)}\right)} + \\
&\hspace{0.4cm}\frac{\widehat{f}(s)}{r_{\text{ext}}}\frac{\Psi_{0,1}\left(r,1,\sqrt{\eta(s)}\right)}{\sqrt{\eta(s)}\Psi_{1,1}\left(1,r_{\text{ext}},\sqrt{\eta(s)}\right)}
\end{split}
\end{eqnarray}
&
\begin{eqnarray}\label{flux_Laplace_caseD}
\begin{split}
\widehat{j}_{2}(s) = \frac{1}{s}
\end{split}
\end{eqnarray}

\\
\bottomrule[0.7pt]
\end{tabular}
\begin{tablenotes}
\item[a] $^*$$\Psi_{m,n}\left(\Phi,\Upsilon,x\right) = K_{m}(\Phi x)I_{n}(\Upsilon x) + (-1)^{m+n+1}I_{m}(\Phi x)K_{n}(\Upsilon x)$ and $\eta(s) = \{[s\omega(1-\omega) + \lambda] s\}/[s(1-\omega) + \lambda]$.
\end{tablenotes}
\end{table}

\bibliography{refHankel}
\end{document}